\documentclass[xypic,amssymb,rawfonts,amsthm,amsmath,amsfonts,syntonly,11pt]{amsart}

\usepackage{amssymb}
\usepackage{hyperref}
\usepackage{mathrsfs}
\usepackage[notref,notcite,final]{showkeys}

\newcounter{thecounter}
\numberwithin{thecounter}{section}
\newtheorem{lemma}[thecounter]{Lemma}
\newtheorem{prop}[thecounter]{Proposition}
\newtheorem{thm}[thecounter]{Theorem}

\newtheorem{Th}{Theorem}

\theoremstyle{definition}

\newtheorem{rem}[thecounter]{Remark}

\numberwithin{equation}{section}

\newcommand*\InsertTheoremBreak{
  \begingroup 
    \setlength\itemsep{0pt}
    \setlength\parsep{0pt}
    \item[\vbox{\null}]
  \endgroup
}

\newcommand{\map}{\operatorname{map}}

\newcommand{\aut}{\operatorname{aut}}
\newcommand{\Aut}{\operatorname{Aut}}

\newcommand{\Inn}{\operatorname{Inn}}
\newcommand{\Out}{\operatorname{Out}}
\newcommand{\End}{\operatorname{End}}

\newcommand{\im}{\operatorname{im}}

\newcommand{\id}{\operatorname{id}}
\newcommand{\res}{\operatorname{res}}
\newcommand{\Der}{\operatorname{Der}}
\newcommand{\PDer}{\operatorname{PDer}}

\newcommand{\cC}{{\mathcal{C}}}

\newcommand{\twocom}{\hat{{}_2}}

\newcommand{\GL}{\operatorname{GL}}
\newcommand{\SO}{\operatorname{SO}}

\newcommand{\SU}{\operatorname{SU}}
\newcommand{\Sp}{\operatorname{Sp}}

\newcommand{\DI}{\operatorname{DI}}

\newcommand{\rad}{\operatorname{rad}}

\newcommand{\Spin}{\operatorname{Spin}}

\newcommand{\luft}{\medskip\par\noindent}
\newcommand{\QED}{\hfill {\bf $\square$}\luft}

\newcommand{\Q}{{\mathbb {Q}}}

\newcommand{\F}{{\mathbb {F}}}
\newcommand{\Z}{{\mathbb {Z}}}

\newcommand{\D}{{\mathbf D}}

\newcommand{\N}{{\mathcal{N}}}

\newcommand{\C}{{\mathbf{C}}}
\newcommand{\cZ}{{\mathcal{Z}}}

\newcommand{\beq}{\begin{eqnarray*}}
\newcommand{\eeq}{\end{eqnarray*}}

\newcommand{\tuborg}{\left\{\begin{array}{ll}}
\newcommand{\sluttuborg}{\end{array}\right.}

\newcommand{\U}{\mathrm{U}}

\input xypic
\xyoption{all}
\CompileMatrices       

\newcommand{\dN}{\breve \N}
\newcommand{\dZ}{\breve \cZ}
\newcommand{\dT}{\breve T}
\newcommand{\dnu}{\breve \nu}
\newcommand{\daut}{\breve{\text{aut}}}

\renewcommand{\phi}{\varphi}

\begin{document}

\title{Automorphisms of $p$-compact groups and their root data}


\author[K. Andersen]{Kasper K. S. Andersen}
\author[J. Grodal]{Jesper Grodal}

\thanks{The second named author was partially supported by NSF grant DMS-0354633 and the Alfred P. Sloan foundation}


\address{Department of Mathematical Sciences, University of Aarhus}

\email{kksa@imf.au.dk}

\address{Department of Mathematics, University of Chicago, Chicago, IL 60637, USA}
\email{jg@math.uchicago.edu}

\begin{abstract}
We construct a model for the space of automorphisms of a connected
$p$-compact group in terms of the space of automorphisms of its maximal
torus normalizer and its root datum.  As a consequence we show that any
homomorphism to the outer automorphism group of a $p$-compact group
can be lifted to a group action, analogous to a classical theorem of de
Siebenthal for compact Lie groups. The model of this paper is used in
a crucial way in our paper ``The classification of
$2$-compact groups'', where we prove the conjectured classification of
$2$-compact groups and determine their automorphism spaces.
\end{abstract}

\maketitle


\section{Introduction}

A $p$-compact group is a homotopy theoretic version of a compact Lie
group at a prime $p$, and consists of a triple $(BX,X,e:X
\xrightarrow{\simeq} \Omega BX)$, where $BX$ is a pointed $p$-complete
space, $X$ is a space with finite $\F_p$-cohomology, and $e$ is a
homotopy equivalence between $X$ and the based loop space on $BX$
\cite{DW94,dwyer98,AGMV03}. The goal of this paper is to give an
algebraic model, which we denote $B\aut(\D_X)$, of $B\Aut(BX)$, the
classifying space of the topological monoid $\Aut(BX)$ of
self-homotopy equivalences of $BX$, for a connected $p$-compact group
$X$. We construct a natural map $\Phi: B\Aut(BX) \to B\aut(\D_X)$ and
prove here that $\Phi$ induces an isomorphism on $\pi_i$, $i>1$
(Theorem~\ref{modelthm}). It is proved to be an isomorphism also on
$\pi_1$ in our sequel paper \cite{AG05classification}, where we prove
the conjectured classification of $2$-compact groups using this model. As a
consequence of the properties of $\Phi$ proved in this paper, we
deduce the following result, which can be seen as a generalization of
a classical theorem of de Siebenthal \cite[Ch.~I, \S 2, no.\
2]{desiebenthal56} \cite[\S 4, no.~10]{bo9} for compact connected Lie
groups (see also Theorem~\ref{plingeling}).

Let $\cZ(X)$ denote the center of the $p$-compact group $X$, $\dZ(X)$
its discrete approximation \cite{dw:center}, and $\Out(BX)$ the group
of components of $\Aut(BX)$.

\begin{Th} \label{kinvthm}
For any connected $p$-compact group $X$, the fibration 
$$
B^2\cZ(X) \to B\Aut(BX) \to B\Out(BX)
$$
has a section. In particular, for a discrete group $\Gamma$,
any group homomorphism $\Gamma \to \Out(BX)$ lifts
to a group action $B\Gamma \to B\Aut(BX)$, and, if
$\widetilde{H}^*(\Gamma;\Q) =0$, the set of liftings is a non-trivial
$H^2(\Gamma;\dZ(X))$-torsor.
\end{Th}

The theorem can fail without the assumption that $X$ is connected, the
smallest examples being $X=D_{16}$ or $Q_{16}$, the dihedral and
quaternion groups of order $16$. Like for Lie groups, the theorem can
be used to construct group actions and hence to construct new
$p$-compact groups or $p$-local finite groups from old, using homotopy
fixed point methods; see e.g., \cite{DW94, BM05}.

\smallskip

Constructing the model requires understanding the relationship between the {\em
  automorphisms} of the maximal torus normalizer $\N_X$ and of the
root datum $\D_X$. On the level of {\em objects} $\D_X$ can be
constructed from $\N_X$ and vice versa: For compact Lie groups this
was proved in an influential 1966 paper by Tits \cite{tits66} and it
was recently generalized to $2$-compact groups by Dwyer and Wilkerson
\cite{DW05}; for $p$-compact groups, $p$ odd, it follows trivially
from a theorem of the first named author \cite[Thm.~1.2]{kksa:thesis}.
On the level of automorphisms the problem is that the outer automorphism groups of
$\N_X$ and $\D_X$ in general differ, which means that we have to
introduce an extra ingredient: In Theorem~\ref{mainthm2cg} we prove
that the outer automorphism group of $\D_X$ naturally corresponds to a
subgroup of the outer automorphism group of $\N_X$, which permutes certain
``root subgroups''. This is the main technical result of the
paper. Even in the case of compact Lie groups the result (stated as
Theorem~\ref{mainthmlie}) appears to be new, though many of the
ingredients are hidden in Tits' paper mentioned above.

\smallskip

We now embark on describing our results more precisely, for which we
need some notation.
For a principal ideal domain $R$, an $R${\em -root datum} $\D$ is defined to
be a triple $(W,L,\{Rb_\sigma\})$, where $L$ is a free
$R$-module of finite rank, $W \subseteq \Aut_R(L)$ is a finite
subgroup generated by reflections (i.e., elements $\sigma$ such that
$1-\sigma \in \End_{R}(L)$ has rank one), and $\{Rb_\sigma\}$
is a collection of rank one submodules of $L$, indexed by the
reflections $\sigma$ in $W$, satisfying
$$
\im(1-\sigma) \subseteq Rb_\sigma
\mbox{ and } w (Rb_\sigma) = Rb_{w\sigma w^{-1}} \mbox{ for all } w \in W.
$$
The element $b_\sigma \in L$, called the {\em coroot} corresponding to
$\sigma$, is determined up to a unit in $R$. Together with $\sigma$ it
determines a {\em root} $\beta_\sigma: L \to R$ via the formula
$$
\sigma(x) = x + \beta_\sigma(x)b_\sigma.
$$
For $R = \Z$ there is a one-to-one correspondence between $\Z$-root
data and classically defined root data, by associating
$(L, \{\pm b_\sigma \}, L^*, \{\pm \beta_\sigma \})$ to
$(W,L,\{\Z b_\sigma\})$; see \cite[Prop.~2.16]{DW05}.
It is easy to see that we always have $Rb_\sigma \subseteq
\ker(1+\sigma + \dots + \sigma^{|\sigma|-1}: L \to L)$, so given $\sigma$ the
possibilities for $Rb_\sigma$ are in bijection with the
cyclic $R$-submodules of $H^1(\left< \sigma \right>;L)$. In particular
if $H^1(\left< \sigma \right>;L) = 0$ for all reflections $\sigma$, as
is the case for $R = \Z_p$, $p$ odd, the root datum $\D = (W, L,
\{Rb_\sigma\})$ is uniquely determined by the underlying reflection group $(W,L)$. 
In fact a coup d'{\oe}il at the classification of finite $\Z$- or $\Z_2$-reflection
groups \cite[Thms.~11.1 and 11.5]{AGMV03} reveals that also for $R =
\Z$ or $\Z_2$ the root datum is determined by $(W,L)$, {\em unless}
$(W,L)$ contains direct factors isomorphic to
$(W_{\SO(2n+1)},L_{\SO(2n+1)})  \cong (W_{\Sp(n)},L_{\Sp(n)})$---these
exceptions, however,
cannot easily be ignored since $\SO(3)$ and $\SU(2)$ are
ubiquitous. Lastly we remark that e.g., for $R
= \Z$ or $\Z_p$ one can instead of the collection $\{Rb_\sigma\}$
equivalently consider their span, the {\em coroot lattice}, $L_0 =
+_\sigma Rb_\sigma \subseteq L$, as is easily seen; this was the definition given in
\cite[\S 1]{AGMV03}, under the name ``$R$-reflection datum''.

Two $R$-root data $\D = (W,L,\{Rb_\sigma\})$ and $\D' =
(W',L',\{Rb'_{\sigma'}\})$ are said to be isomorphic if there exists an
isomorphism $\phi: L \to L'$ such that $\phi W\phi^{-1} = W'$ and
$\phi(Rb_\sigma) = Rb'_{\phi \sigma \phi^{-1}}$. In particular the
{\em automorphism group} of $\D$ is given by $\Aut(\D) = \{ \phi \in
N_{\Aut_R(L)}(W) | \phi(Rb_\sigma) = Rb_{\phi \sigma \phi^{-1}} \}$ and we
define the {\em outer automorphism group} as $\Out(\D) = \Aut(\D)/W$.

Tits constructed in \cite{tits66} for any
$\Z$-root datum $\D$ an extension 
$$
1 \to T \to \nu(\D) \to W \to 1,
$$
where $T = L \otimes_{\Z} S^1$, such that $N_G(T) \cong \nu(\D_G)$ for any
compact connected Lie group $G$ with root datum $\D_G$
(see Theorem~\ref{sigmaident}\eqref{maplie}). The group $\nu(\D)$ will
be called the {\em maximal torus normalizer} associated to
$\D$. Conversely the $\Z$-root
datum $\D$ can be recovered from $\nu(\D)$, see
Proposition~\ref{reconstr}. For each reflection $\sigma \in W$ one can
algebraically construct a canonical rank one ``root subgroup''
$N_\sigma = \nu(\D)_\sigma$ inside $N = \nu(\D)$, and similarly there
are root subgroups $N_G(T)_{\sigma}$ of $N_G(T)$. A more precise
version of Tits' result is that there is an essentially  {\em unique}
isomorphism between $N_G(T)$ and $\nu(\D)$ sending $N_G(T)_\sigma$ to
$\nu(\D)_\sigma$ (see Theorem~\ref{sigmaident}\eqref{maplie}).
Define $\Aut(N,\{N_\sigma\}) = \{
\phi \in \Aut(N) \, | \, \phi(N_\sigma) = N_{\phi \sigma
  \phi^{-1}}\}$, and $\Out(N,\{N_\sigma\}) =
\Aut(N,\{N_\sigma\})/\Inn(N)$ (where $\Aut$ and $\Out$ here have their
usual meaning for Lie groups!). As usual $H^1(W;T) =
\Der(W,T)/\PDer(W,T)$, derivations modulo principal derivations, and
$H^1(W;T)$ naturally embeds in $\Out(N)$ via the homomorphism $[f]\mapsto [\phi_f]$,
where $\phi_f\in \Aut(N)$ is given by $\phi_f(x) = f(\bar x)x$.
The following is our main result in the Lie group case.

\begin{Th} \label{mainthmlie}
Let $\D = (W,L,\{\Z b_\sigma\})$ be a $\Z$-root datum with corresponding
maximal torus normalizer $N = \nu(\D)$. There
is a short exact sequence
$$
1 \to H^1(W;T) \to \Out(N) \to \Out(\D) \to 1,
$$
which has a {\em canonical} splitting $s: \Out(\D) \to \Out(N)$ with image
$\Out(N,\{N_\sigma\})$.

Furthermore, if $G$ is a compact connected Lie group with root datum
$\D_G$, the splitting $s$ fits into the following commutative diagram of isomorphisms
$$ 
\xymatrix{
 & \Out(N,\{N_\sigma\}) \ar@/^/[dr]^{\res}_\cong & \\
\Out(G) \ar[ur]^{\res}_\cong \ar[rr]^{\res}_\cong & & \Out(\D_G).
\ar@/^/[ul]^{s}
}
$$
\end{Th}

A split exact sequence as above was first established by H\"ammerli
\cite[Thm.~1.2]{hammerli02}, using the existence of an ambient Lie
group and case-by-case arguments, but without identifying the image of
the splitting, the key point for our purposes. 

\smallskip

For $p$ odd, define
 the {\em discrete maximal torus normalizer}
$\dnu(\D)$ of a
$\Z_p$-root datum $\D$ to be the semidirect
product $\dT \rtimes W$, where $\dT = L \otimes_{\Z} \Q/\Z$. 
By a theorem of the first named author \cite[Thm.~1.2]{kksa:thesis},
$\dN_X \cong \dnu(\D_X)$, for any connected $p$-compact group $X$,
where $\D_X$ is the unique $\Z_p$-root datum with underlying
reflection group $(W_X, L_X)$.

For $p=2$, Dwyer-Wilkerson recently constructed in
\cite{DW05} a group $\dnu(\D)$ (again called the discrete maximal
torus normalizer) from a $\Z_2$-root datum $\D$, generalizing Tits'
construction, as the middle term in an extension
$$
1 \to \dT \to \dnu(\D) \to W \to 1,
$$
and showed how to associate a $\Z_2$-root datum $\D_X$ to any
connected $2$-compact group $X$, such that $\dN_X \cong
\dnu(\D_X)$. As for $\Z$-root data we construct in Section~\ref{normalext-sec}
algebraic root subgroups $\dN_\sigma = \dnu(\D)_\sigma$ of $\dN =
\dnu(\D)$ for any $\Z_p$-root datum $\D$ and similarly one can define
subgroups $(\dN_X)_\sigma$ of $\dN_X$. We strengthen the result of
Dwyer-Wilkerson in
Section~\ref{normalext-sec} (Theorem~\ref{sigmaident}\eqref{mappcg}),
keeping track of root subgroups, and this sharpening is key to our results.
At the same time we remove an unfortunate reliance in \cite{DW05} on a
classification of $2$-compact groups of semi-simple rank $2$.

As above define $\Aut(\dN,\{\dN_\sigma\}) = \{
\phi \in \Aut(\dN) \, | \, \phi(\dN_\sigma) = \dN_{\phi \sigma
  \phi^{-1}}\}$ and $\Out(\dN,\{\dN_\sigma\}) =
\Aut(\dN,\{\dN_\sigma\})/\Inn(\dN)$. For $p$-compact groups, we prove the following analog
of Theorem~\ref{mainthmlie}. 

\begin{Th} \label{mainthm2cg}
Let $\D = (W,L,\{\Z_p b_\sigma\})$ be a $\Z_p$-root datum with corresponding
discrete maximal torus normalizer $\dN = \dnu(\D)$. There
is a short exact sequence
$$
1 \to H^1(W;\dT) \to \Out(\dN) \to \Out(\D) \to 1,
$$
which has a {\em canonical} splitting $s: \Out(\D) \to \Out(\dN)$ with image
$\Out(\dN,\{\dN_\sigma\})$.

Furthermore, if $X$ is a connected $p$-compact group with discrete
maximal torus normalizer $\dN = \dN_X$, the canonical homomorphism $\Phi:
\Out(BX) \to \Out(\dN)$ has image contained in
$\Out(\dN,\{\dN_\sigma\})$.
\end{Th}

Here $\Phi$ is $\pi_1$ of the ``Adams-Mahmud'' map $\Phi: B\Aut(BX)
\to B\Aut(B\N_X)$ given by lifting self-equivalences of $BX$ to
self-equivalences of $B\N_X$, see \cite[Lem.~4.1 and
Prop.~5.1]{AGMV03}. For $p$ odd, Theorem~\ref{mainthm2cg} degenerates
since in this case $H^1(W;\dT) = 0$ by \cite[Thm.~3.3]{kksa:thesis}.

We now define the space $B\aut(\D)$, which provides an algebraic
model for $B\Aut(BX)$, depending only on $\D = \D_X$. For $p$ odd, set
$B\aut(\D) = B\Aut(B\N_X)$, and recall that in this case 
$\Phi: B\Aut(BX) \to B\aut(\D)$ is a homotopy equivalence by
\cite[Thm.~1.4]{AGMV03}. For $p=2$ we need to modify $B\Aut(B\N_X)$ to
get a model, since the $H^1$-term in Theorem~\ref{mainthm2cg} can be non-zero.
Let $Y$ be the covering space of $B\Aut(B\N_X)$ corresponding to the
subgroup $\Out(\dN_X,\{(\dN_X)_\sigma\})$ of the fundamental
group. The space $B\aut(\D)$ is obtained from $Y$ by modifying
$\pi_2(Y)$ by killing certain $\Z/2$-summands, one for each direct factor of $\D$
isomorphic to $\D_{\SO(2n+1)\twocom}$; we refer to
Section~\ref{kinv-sec1} for a precise description of this. With
$\cZ(\D)$ denoting the center of $\D$, defined in
Section~\ref{kinv-sec1}, we prove in Section~\ref{kinv-sec2}  that we
can identify the homotopy type of $B\aut(\D)$ with $(B^2\cZ(\D))_{h\Out(\D)}$.

\begin{Th} \label{modelthm}
For a connected $p$-compact group $X$ with
root datum $\D_X$, the Adams-Mahmud map
induces a natural map
$$
\Phi: B\Aut(BX) \to (B^2\cZ(\D_X))_{h\Out(\D_X)}
$$
which is an isomorphism on $\pi_i$ for $i>1$ and is the canonical map
$\Phi: \Out(BX) \to \Out(\D_X)$ on $\pi_1$.
\end{Th}

From this Theorem~\ref{kinvthm} follows easily. As mentioned, we prove
in our sequel paper \cite{AG05classification}, as part of our
inductive proof of the classification of $2$-compact groups, that
$\Phi$ also induces an isomorphism on $\pi_1$ for $p=2$. For $p$ odd,
$\Phi$ induces an isomorphism on $\pi_1$ by our work with M{\o}ller
and Viruel \cite[Thm~1.1]{AGMV03}.

\smallskip

{\em Organization of the paper:}
In Section~\ref{reflext-sec} we recall the reflection extension of
Dwyer-Wilkerson, and establish some further properties, which we use
in Section~\ref{normalext-sec} to construct the root subgroups
$\nu(\D)_\sigma$. In Section~\ref{proof-sec} we prove
Theorems~\ref{mainthmlie} and \ref{mainthm2cg},
and in Section~\ref{kinv-sec1} we construct the space $B\aut(\D)$ which we
use in Section~\ref{kinv-sec2} to prove Theorems~\ref{kinvthm} and \ref{modelthm}.

\smallskip

{\em Notation:}
For a $\Z\left<\sigma\right>$-module $A$ we define
$A^-(\sigma) = \ker(\sum_{i=0}^{|\sigma|-1}\sigma^i : A \to A)$ and
$A^+(\sigma) = \ker(1-\sigma: A \to A)$. If $A$ is an abelian
compact Lie group, $A_0$ denotes its identity component. For a compact Lie
group $G$, $\Aut(G)$ denotes the group of Lie group automorphisms of
$G$. If $Y$ is a topological space, $\Aut(Y)$ denotes the topological monoid of
self-homotopy equivalences of $Y$. A $\Z_p$-root datum $\D'$ is said to be of
{\em Coxeter type} if $\D' \cong \D \otimes_\Z \Z_p = (W,L\otimes_{\Z}
\Z_p,\{\Z_p \otimes_\Z \Z b_\sigma\})$ for a $\Z$-root datum $\D =
(W,L,\{\Z b_\sigma\})$. We call a $\Z_p$-root datum $\D = (W, L, \{\Z
b_\sigma\})$ {\em exotic} if $L\otimes_{\Z} \Q$ is an irreducible
$\Q_p[W]$-module and $\D$ is not of Coxeter type.
A {\em subdatum} of a $\Z_p$-root datum $\D = (W, L, \{\Z_p
b_\sigma\})$ is a $\Z_p$-root datum of the form $(W', L, \{\Z_p
b_\sigma\}_{\sigma\in \Sigma'})$ where $(W', L)$ is a reflection
subgroup of $(W, L)$ and $\Sigma'$ is the set of reflections in $W'$.
We will freely use the terminology
of $p$-compact groups, though we try to give concrete references for
the facts we use---we refer the reader to \cite{DW94,dwyer98,AGMV03}
for background information.


\section{The reflection extension} \label{reflext-sec}

In this section we introduce
certain subextensions $\rho(W)_\sigma$ of the {\em reflection
  extension} $\rho(W)$ of Dwyer-Wilkerson \cite{DW05} (recalled below)
for any reflection $\sigma$ in a finite $\Q_2$-reflection group $W$,
and describe their basic properties (Proposition~\ref{gen210}). The
results generalize work of Tits \cite{tits66} over $\Q$, and are key
to constructing the root subgroups $\nu(\D)_\sigma$ of the normalizer
extension $\nu(\D)$ in the next section.

We begin with some definitions and recollections. Let $W$ be a finite
$\Q_2$-reflection group and let $\Sigma = \bigcup_i \Sigma_i$ denote
the partition of the set of reflections $\Sigma$ in $W$ into conjugacy
classes. Choose for each conjugacy class $\Sigma_i$ an element
$\tau_i\in \Sigma_i$. Since $-1$ is the only non-trivial element of
finite order in $\Q_2^\times$ we may find an element $a_i\neq 0$ in
the $\Q_2$-vector space underlying $W$ such that $\tau_i(a_i) =
-a_i$. Let $C_i = C_W(\tau_i)$ and $C_i^\perp$ be the stabilizer of
$a_i$ in $W$. It is clear that we get a direct product decomposition
$C_i = \left<\tau_i\right> \times C_i^\perp \cong \Z/2 \times
C_i^\perp$. Now consider the extension
$$
1 \to \Z \xrightarrow{(\cdot 2,\cdot 0)} \Z \times C_i^\perp \to \Z/2
\times C_i^\perp \to 1,
$$
where $\cdot 0$ denotes the trivial homomorphism.
This extension produces an element in $H^2(C_i;\Z)$. By Shapiro's
lemma we have $H^2(C_i;\Z) \cong H^2(W;\Z[\Sigma_i])$ and hence we get
an extension $\rho_i(W)$ of $W$ by $\Z[\Sigma_i]$. Following
Dwyer-Wilkerson \cite{DW05} we define the {\em reflection extension}
$\rho(W)$ of $W$ as the sum of the extensions $\rho_i(W)$. Noting that
$\bigoplus_i \Z[\Sigma_i] \cong \Z[\Sigma]$, we see that the
reflection extension has the form $1 \to \Z[\Sigma] \to \rho(W)
\xrightarrow{\pi} W \to 1$. It is clear that the extensions $\rho_i(W)$ are
well-defined up to equivalence of extensions and hence $\rho(W)$ is
also well-defined.

\begin{lemma} \label{products}
Let $W_1$ and $W_2$ be finite $\Q_2$-reflection groups and set $W =
W_1 \times W_2$. Then the reflection extension of $W$ splits as a
product
$$
\rho(W) \cong \rho(W_1) \times \rho(W_2).
$$
\end{lemma}

\begin{proof}
Let $\tau_i$ be a reflection in $W_1$. For
the construction of $\rho(W_1)$ we have to consider the extension
$$
1 \to \Z \xrightarrow{(\cdot 2,\cdot 0)} \Z \times C_i^\perp \to \Z/2
\times C_i^\perp \to 1,
$$
where $C_i = C_{W_1}(\tau_i)$. Since $C_W(\tau_i) = C_i \times W_2$
and $C_W(\tau_i)^\perp = C_i^\perp \times W_2$ it is clear that the
associated sequence to be used in the construction of $\rho(W)$ equals
$$
1 \to \Z \xrightarrow{(\cdot 2,\cdot 0)} \Z \times (C_i^\perp \times W_2) \to \Z/2
\times (C_i^\perp \times W_2) \to 1.
$$
Via Shapiro's lemma the first extension corresponds to an extension
$$
1 \to \Z[\Sigma_i] \to \rho_i(W_1) \to W_1 \to 1.
$$
It is now clear that the associated sequence for $W$ must be
$$
1 \to \Z[\Sigma_i] \to \rho_i(W_1) \times W_2 \to W_1 \times W_2 \to 1.
$$
From this it follows that the reflection extension for $W$ must be the
sum of the extensions $\rho(W_1) \times W_2$ and $W_1 \times
\rho(W_2)$ which equals the extension $\rho(W_1) \times \rho(W_2)$ as claimed.
\end{proof}

We also need the following lemma, which can be extracted from \cite{DW05}.

\begin{lemma}[{\cite[Pf.~of~Lem.~10.1]{DW05}}] \label{pullback}
Let $W$ be a finite $\Q_2$-reflection group, and let $W_1$ be a
reflection subgroup. Let $\Sigma$ be the set of reflections in $W$ and
$\Sigma_1$ the set of reflections in $W_1$. Then the pullback of the
reflection extension
$1 \to \Z[\Sigma] \to \rho(W) \to W \to 1$ along the inclusion $W_1 \to
W$ is the sum of the reflection extension $1 \to \Z[\Sigma_1] \to
\rho(W_1) \to W_1 \to 1$ and the semidirect product
$1 \to \Z[\Sigma\setminus \Sigma_1] \to \Z[\Sigma\setminus \Sigma_1]
\rtimes W_1 \to W_1 \to 1$.\QED
\end{lemma}

We can now prove the main statement of this section, which generalizes
\cite[Prop.~2.10]{tits66} to finite $\Q_2$-reflection groups. 

\begin{prop} \label{gen210}
Let $W$ be a finite $\Q_2$-reflection group and $\rho(W)$ the
associated reflection extension. For a reflection $\sigma\in W$ we
define $C_\sigma = (1-\sigma) \Z[\Sigma]$, $Q_\sigma = \{ x\in
\rho(W)\, |\, x^2 = \sigma\in \Z[\Sigma] \}$,  and $\rho(W)_\sigma =
\left< Q_\sigma \right>$. Then
\begin{enumerate}
\item \label{conj}
For $x\in \rho(W)$ with $w=\pi(x)$ we have $x C_\sigma x^{-1} =
C_{w\sigma w^{-1}}$ and $x Q_\sigma x^{-1} = Q_{w\sigma w^{-1}}$.
\item \label{invol}
$\rho(W)$ contains no involutions.
\item \label{proj}
$\pi(Q_\sigma) = \{\sigma\}$.
\item \label{Csigma}
$C_\sigma$ is a subgroup of $\Z[\Sigma]$ and we have $C_\sigma =
\ker(\Z[\Sigma] \xrightarrow{1+\sigma} \Z[\Sigma])$.
\item \label{coset}
$Q_\sigma$ is a coset of $C_\sigma$ in $\rho(W)$; more precisely, for
$x\in Q_\sigma$ we have $Q_\sigma = x C_\sigma = C_\sigma x$.
\item \label{ses}
We have a short exact sequence $1 \to C_\sigma \oplus \Z \sigma \to
\rho(W)_\sigma \xrightarrow{\pi} \left<\sigma\right> \to 1$.
\end{enumerate}
\end{prop}

\begin{proof}
Part \eqref{conj} is obvious. By Lemma~\ref{products} we see that if $W =
W_1\times W_2$ is a product of two finite $\Q_2$-reflection groups,
then $\rho(W) = \rho(W_1) \times \rho(W_2)$. Thus $W$ satisfies
\eqref{invol} if and only if $W_1$ and $W_2$ does.

Moreover, if $\sigma\in W_1$ is a reflection we get
$$
Q_\sigma(W) = Q_\sigma(W_1) \times \{x\in \rho(W_2)\, |\, x^2=1 \},
$$
where $Q_\sigma(W)$ and $Q_\sigma(W_1)$ denotes $Q_\sigma$ defined with
respect to $W$ and $W_1$ respectively. In particular $W$ satisfies
\eqref{proj} for a fixed $\sigma\in W_1$ if and only if $W_1$ satisfies
\eqref{proj} for this $\sigma$ and $W_2$ satisfies \eqref{invol}.

Since any finite $\Q_2$-reflection group is a product of a Coxeter
group and a number of copies of $W_{\DI(4)}$ (cf.\ \cite{CE74}), the
above remarks shows that it suffices to prove \eqref{invol} and
\eqref{proj} in the Coxeter case and in the case of $W_{\DI(4)}$.

We first deal with the Coxeter case. Here \eqref{proj} holds by
\cite[Prop.~2.10(b)]{tits66}. Now let $W$ be any Coxeter group and let
$W'$ be a non-trivial Coxeter group. Choose any reflection $\sigma\in
W'$. Since \eqref{proj} holds for the Coxeter group $W' \times W$ with
respect to $\sigma$ we obtain \eqref{invol} for $W$ by the above. This
proves \eqref{invol} and \eqref{proj} in the Coxeter case.

In the case $W=W_{\DI(4)}$, there is a reflection subgroup $W_1$ of $W$
isomorphic to $W_{\Spin(7)}$ (cf.\ e.g.\ \cite[Pf.~of~Prop.~9.12, $\DI(4)$
case]{DW05}). It is easily checked that any element $w\in W$ with $w^2
= 1$ is conjugate to an element in $W_1$. Assume first that $x\in
\rho(W)$ with $x^2 = 1$. Then $\pi(x) \in W$ satisfies $\pi(x)^2 = 1$ so
that up to a conjugation in $\rho(W)$ we may assume $\pi(x) \in
W_1$, i.e., that $x$ belongs to the pullback of $\rho(W)$ along the
inclusion $W_1 \to W$. Under the isomorphism given by
Lemma~\ref{pullback}, the element $x$ corresponds to a pair of
elements $(x_1,x_2)$ with $x_1\in \rho(W_1)$ and $x_2 \in
\Z[\Sigma\setminus \Sigma_1] \rtimes W_1$ having the same image under
the projections to $W_1$. Since $x^2=1$ we have $x_1^2=1$ and
$x_2^2=1$. Since we have already established \eqref{invol} for the
Coxeter group $W_1$ we conclude that $x_1=1$. Hence $x_2$ projects
to the identity in $W_1$, so $x_2$ is contained in the subgroup
$\Z[\Sigma\setminus \Sigma_1]$ of $\Z[\Sigma\setminus \Sigma_1] \rtimes
W_1$. Since this subgroup contains no involutions we see that
$x_2=1$. Therefore $x=1$ which proves \eqref{invol} for $W_{\DI(4)}$.

To prove \eqref{proj} for $W_{\DI(4)}$, note first that by
Lemma~\ref{pullback} $\rho(W_1)$ may be considered as a
subgroup of the pullback of $\rho(W)$ along $W_1 \to W$ since this
pullback is the sum of the reflection extension of $W_1$ and a
semidirect product. Thus $\rho(W_1)$ may be viewed as a subgroup of
$\rho(W)$. Since \eqref{proj} holds for the Coxeter group $W_1$ we get
$\sigma \in \pi(Q_\sigma)$ for $\sigma\in \Sigma_1$. Hence $\sigma \in
\pi(Q_\sigma)$ for all $\sigma\in \Sigma$ by \eqref{conj} as any
reflection in $W$ is conjugate to one in $W_1$.

Conversely, assume that $x\in Q_\sigma$ for some $\sigma\in
\Sigma$. As above $\pi(x) \in W$ satisfies $\pi(x)^2 = 1$, so up
to conjugation in $\rho(W)$ we may assume that $\pi(x) \in W_1$, i.e.\ that
$x$ belongs to the pullback of $\rho(W)$ along the inclusion $W_1 \to
W$. As above $x$ corresponds to a pair of elements $(x_1,x_2)$ with
$x_1\in \rho(W_1)$ and $x_2 \in \Z[\Sigma\setminus \Sigma_1] \rtimes
W_1$ having the same image under the projections to $W_1$. Then $x^2 =
\sigma \in \Z[\Sigma]$ corresponds to $(x_1^2,x_2^2)$. Hence if $\sigma\in
\Sigma_1$ we have $x_1^2 = \sigma$ (and $x_2^2=1$) and if
$\sigma\notin \Sigma_1$ we have $x_2^2=\sigma$ (and
$x_1^2=1$). However the second case cannot occur since no element in
$\Z[\Sigma \setminus \Sigma_1]\rtimes W_1$ has square $\sigma$ (the
square of an element $(x,w)$ in the semidirect product equals
$(x+w\cdot x,w^2)$ and the image of $x+w\cdot x$ under the
augmentation $\Z[\Sigma \setminus \Sigma_1] \to \Z$ is even).
We thus conclude that $\sigma\in \Sigma_1$ and $x_1^2 = \sigma \in
\Sigma_1$. Thus $x_1$ belongs to the subgroup $Q_\sigma(W_1)$ and
hence $x_1$ projects to $\sigma\in W_1$ since \eqref{proj} holds for
$W_1$. This shows that $\pi(x) = \sigma$ as desired. This proves
\eqref{proj} for $W_{\DI(4)}$ and hence \eqref{invol} and \eqref{proj}
holds in general.

Part \eqref{Csigma} amounts to showing that
$H^1(\left<\sigma\right>;\Z[\Sigma]) = 0$. Since $\sigma$ is an
involution the $\left<\sigma\right>$-module $\Z[\Sigma]$ splits as a
direct sum of $\Z$'s with trivial action and $\Z^2$'s with the action given by
permuting the generators. Since $H^1(\left<\sigma\right>;-) = 0$ in
both cases, the claim follows.

To prove \eqref{coset}, consider $x\in Q_\sigma$. By
\eqref{proj} any element in $Q_\sigma$ has the form $xu$ for
some $u\in \Z[\Sigma]$. The computation $(xu)^2 = x^2 x^{-1}ux u =
\sigma x^{-1}ux u$ shows that for $u\in \Z[\Sigma]$ we have $xu\in
Q_\sigma$ if and only if $x^{-1}ux u=1$. Since $\pi(x)=\sigma$ this is
equivalent to $(1+\sigma)\cdot u = 0$. By \eqref{Csigma} this proves
that $Q_\sigma = x C_\sigma$ and one gets $Q_\sigma = C_\sigma x$
analogously.

Finally, to see \eqref{ses}, note that $\rho(W)_\sigma \xrightarrow{\pi}
\left<\sigma\right>$ is surjective by \eqref{proj}. Fixing $x\in
Q_\sigma$ we have $\rho(W)_\sigma = \left<x,C_\sigma\right>$ by
\eqref{coset}. Since $x^k C_\sigma = C_\sigma x^k$ for $k\in\Z$
by \eqref{coset}, any element in $\rho(W)_\sigma$ has the form $x^k
c$ for some $k\in\Z$ and some $c\in C_\sigma$. As $x\notin
\ker(\pi)$ and $x^2=\sigma \in\ker(\pi)$ it follows that the kernel of
$\rho(W)_\sigma \xrightarrow{\pi} \left<\sigma\right>$ equals
$\left<C_\sigma,\sigma\right>$ proving the proposition since
$\sigma\notin C_\sigma$.
\end{proof}


\section{The normalizer extension} \label{normalext-sec}

In this section we define the normalizer extension and the associated
root subgroups algebraically, using the results and notation of the
previous section. We use this to give a strengthened version, needed
for our purposes, of the result of Dwyer-Wilkerson
\cite[Prop.~1.10]{DW05} on the maximal torus normalizer in a connected
$2$-compact group, where we also keep track of certain root subgroups
defined below. This version
is on a par with the corresponding result for compact Lie
groups by Tits \cite{tits66}.
We furthermore circumvent the use of a classification of connected
$2$-compact groups of semi-simple rank $2$ used in
\cite[Prop.~9.15]{DW05}. (The semi-simple rank of a $p$-compact group
is the rank of its universal cover.) For this
Dwyer-Wilkerson refer to \cite[Thm.~6.1]{BKNP03} which collects a
range of disparate results, including old classification results on
finite $H$-spaces. We prefer to avoid this reliance, since we use the
results of this paper in our general classification of $2$-compact
groups \cite{AG05classification}, and we are able to do this using  a
few low-dimensional group cohomology computations instead.

We start by giving an alternative definition of $\Z$- and $\Z_p$-root
data in terms of {\em markings}. For a $\Z$-root datum $\D = (W, L,
\{\Z b_\sigma\})$, the {\em marking} associated to the reflection
$\sigma\in W$ is the element $h_\sigma = b_\sigma/2 \in T = L
\otimes_\Z S^1$. Conversely, $\Z b_\sigma = \frac{2}{|h_\sigma|}
\ker(L\xrightarrow{1+\sigma} L)$, so one might as well define a
$\Z$-root datum in terms of the $h_\sigma$ instead of the $\Z
b_\sigma$, and we will use these two viewpoints interchangeably
without further comment. (The conditions on $h_\sigma$ corresponding
to the conditions on $\Z b_\sigma$ say that $h_\sigma \in
T_0^-(\sigma)$, $h_\sigma^2=1$ and $h_\sigma\neq 1$ if $\sigma$ acts
non-trivially on ${}_2 T$, cf.\ \cite[Def.~2.12]{DW05}).
This definition of the markings $h_\sigma$ carries over verbatim to
$\Z_p$-root data replacing $T$ by $\dT = L \otimes_{\Z} \Q/\Z$; see \cite[\S
6]{DW05} for the case $p=2$, for $p>2$ we simply have
$h_\sigma = 1\in \dT$.

Next we recall the construction of root data for compact Lie groups
and $p$-compact groups: For a compact connected Lie group $G$ with
maximal torus $T$, maximal torus normalizer $N=N_G(T)$, Weyl group
$W=N/T$ and integral lattice $L=\pi_1(T)$, we define the
{\em topological root subgroup} as
$$
N_\sigma = \{ x \in C_N(T_0^+(\sigma)) \,|\, x \mbox{ is conjugate in
} C_G(T_0^+(\sigma)) \mbox{ to some } y \in T^-_0(\sigma) \},
$$
for any reflection $\sigma\in W$. This is a subgroup of $N$, $\{x^2
\,|\, x \in N_\sigma \setminus T_0^-(\sigma)\}$ consists of a single
element $h_\sigma$ (see \cite[Lem.~5.2]{DW05}), and $\D_G =
(W,L,\{h_\sigma\})$ is a $\Z$-root datum, the root datum of $G$ (see
\cite[Prop.~5.5]{DW05}). One can check that $N_\sigma =
C_N(T_0^+(\sigma)) \cap C_G(T_0^+(\sigma))'$ (the prime denoting
the commutator subgroup), the definition from
\cite[4.1]{tits66}, and that $C_N(T_0^+(\sigma))$ is the pullback of
$N$ along $\left<\sigma\right> \to W$. Note that by definition, the
topological root subgroups are well behaved under conjugation: For
$x\in N$ with image $w\in W$ we have $x N_\sigma x^{-1} = N_{w\sigma
  w^{-1}}$ (cf.\ \cite[Prop.~5.5]{DW05}).

For a connected $2$-compact group $X$, we can use the above definition verbatim
(replacing $G$ by $X$, $N$ by $\dN_X$ and $T$ by $\dT$) to define
{\em topological root subgroups} $(\dN_X)_\sigma$ of $X$ and a $\Z_2$-root
datum $\D_X$, the {\em root datum of $X$}. The verification that these
have the desired properties is given in \cite[\S 9]{DW05} and uses
that any connected $2$-compact group of semi-simple rank $1$ is either
$(\SU(2) \times (S^1)^n)\twocom$, $(\SO(3) \times (S^1)^{n})\twocom$,
or $(\U(2) \times (S^1)^{n})\twocom$ for some $n\geq 0$; this is
essentially classical, and follows easily from \cite{DMW86}.
As above we also have $x \dN_\sigma x^{-1} = \dN_{w\sigma
  w^{-1}}$ for $x\in \dN_X$ with image $w\in W_X$ \cite[Prop.~9.10]{DW05}.

For a connected $p$-compact group $X$, $p$ odd, the first-named author
proved that the extension $1 \to \dT \to \dN_X \to W \to 1$ is split
\cite[Thm.~1.2]{kksa:thesis}, and the splitting is unique up to
conjugation by an element in $\dT$ \cite[Thm.~3.3]{kksa:thesis}. Hence
we can define the {\em topological root subgroups} $(\dN_X)_\sigma$ of
$\dN_X$ as the images of the subgroups $\dT^-(\sigma) \rtimes
\left<\sigma\right> \subseteq \dT \rtimes W$ under the isomorphism
$\dN_X \cong \dT \rtimes W$ of extensions; this is independent of the
chosen isomorphism since two such isomorphisms differ by conjugation by an
element in $\dT$. Note also that $\dT^-(\sigma) = \dT^-_0(\sigma)$
since $H^1(\left<\sigma\right>;\dT) = 0$ for $p$ odd. We also define
the root datum $\D_X$ of $X$ as the unique $\Z_p$-root datum whose
underlying reflection group is $(W_X, L_X)$.

\smallskip

We define, following Dwyer-Wilkerson \cite[\S 3]{DW05}, the
{\em normalizer extension}
$$
1 \to T \to \nu(\D) \to W \to 1
$$ 
of a $\Z$-root datum $\D = (W,L,\{\Z b_\sigma\})$ as the pushforward
of the reflection extension $1 \to \Z[\Sigma] \to \rho(W) \to W \to 1$
constructed in Section~\ref{reflext-sec}, along the $W$-map
$f:\Z[\Sigma] \to T$ given by $\sigma \mapsto h_\sigma$.

For a reflection $\sigma\in W$, define the {\em algebraic root
  subgroup} $\nu(\D)_\sigma$ of $\nu(\D)$ via the extension
$$
1 \to T_0^-(\sigma) \to \nu(\D)_\sigma \to \left<\sigma\right> \to 1
$$
obtained as the pushforward along $f_{|C_\sigma \oplus \Z\sigma}:
C_\sigma \oplus \Z\sigma \to T_0^-(\sigma)$ of the extension $1 \to C_\sigma
\oplus \Z \sigma \to \rho(W)_\sigma \xrightarrow{\pi}
\left<\sigma\right> \to 1$, constructed in Proposition~\ref{gen210}\eqref{ses}.
Note that this makes sense, since $f(C_\sigma \oplus \Z\sigma)$
is contained in $T_0^-(\sigma)$ because $f(C_\sigma) = f(
(1-\sigma)\Z[\Sigma]) \subseteq (1-\sigma) T = T_0^-(\sigma)$
and $f(\sigma) = h_\sigma \in T_0^-(\sigma)$. Also note that by
Proposition~\ref{gen210}\eqref{conj} we have $x \nu(\D)_\sigma x^{-1}
= \nu(\D)_{w\sigma w^{-1}}$ for $x\in \nu(\D)$ with image $w\in W$.

Similarly, for a $\Z_2$-root datum $\D = (W,L,\{\Z_2 b_\sigma\})$ the
{\em normalizer extension} $1\to \dT  \to \dnu(\D) \to W \to 1$ is
defined to be the pushforward of $\rho(W)$ along $f:\Z[\Sigma] \to
\dT$ defined by $\sigma \mapsto h_\sigma$, and again the {\em
  algebraic root subgroup} $\dnu(\D)_\sigma$ is defined via the extension
$1 \to \dT_0^-(\sigma) \to \dnu(\D)_\sigma \to
\left<\sigma\right> \to 1$ obtained as the pushforward of $1 \to C_\sigma
\oplus \Z \sigma \to \rho(W)_\sigma \xrightarrow{\pi}
\left<\sigma\right> \to 1$ along $f_{|C_\sigma \oplus \Z\sigma}
:C_\sigma \oplus \Z\sigma \to \dT_0^-(\sigma)$.

Recall that for a $\Z_p$-root datum $\D = (W,L,\{\Z_p b_\sigma\})$, $p$
odd, we have defined the discrete maximal torus normalizer $\dnu(\D)$ as
$\dT \rtimes W$. Similarly we define the {\em algebraic root
  subgroups} by $\dnu(\D)_\sigma = \dT^-(\sigma) \rtimes
\left<\sigma\right>$.

\begin{rem}
We remark that the name root subgroup for $N_\sigma$ perhaps is a bit
unfortunate: If $G(\C)$ is a connected reductive algebraic group over
$\C$ and $G$ is the corresponding maximal compact subgroup, then
$N_\sigma$ with respect to $G$ is the maximal compact subgroup in the
normalizer of $T(\C)^-_0(\sigma)$ in $\left< U_\alpha,
  U_{-\alpha}\right>$, where $\pm \alpha$ are the two roots associated
to $\sigma$ and $U_\alpha$ is what is usually called the root subgroup
corresponding to $\alpha$ in the algebraic group $G(\C)$. (See
\cite[Thm.~IV.13.18(4)(d)]{borel91},
\cite[Prop.~8.1.1(i)]{spr:linalggrp} or
\cite[Thm.~26.3]{humphreys75}.)
\end{rem}

\smallskip

We can now state an improved version of the main theorem of
 Dwyer-Wilkerson \cite{DW05}, which says that the
algebraic and topological definitions of the maximal torus normalizer
and the root subgroups coincide.

\begin{thm} \label{sigmaident}
We have the following identifications.
\begin{enumerate}
\item \label{maplie} \cite[Thm.~4.4]{tits66}
If $G$ is a compact connected Lie group with root datum $\D_G =
(W,L,\{\Z b_\sigma\})$, then there exists an isomorphism of extensions
$$
\xymatrix{
1 \ar[r] & T \ar@{=}[d] \ar[r]& \nu(\D_G) \ar[d]^\cong \ar[r] & W
\ar@{=}[d] \ar[r]& 1\\
1 \ar[r] & T \ar[r]& N_G(T) \ar[r]& W \ar[r]& 1,
}
$$
taking $\nu(\D_G)_\sigma$ to $N_G(T)_\sigma$ for all $\sigma$, and
this specifies the isomorphism uniquely up to conjugation by an
element in $T$.
\item \label{mappcg} (compare \cite[Prop.~1.10]{DW05})
If $X$ is a connected $p$-compact group with root datum $\D_X =
(W,L,\{\Z_p b_\sigma\})$ and discrete maximal torus normalizer $\dN_X$,
then the result of \eqref{maplie} continues to hold replacing $G$ by
$X$, $T$ by $\dT$, $\nu(\D_G)$ by $\dnu(\D_X)$ and $N_G(T)$ by $\dN_X$.
\end{enumerate}
\end{thm}

Before the proof we need to establish some properties of the algebraic
root subgroups. The first lemma tells us how to recover the elements
$h_\sigma$ from the root subgroups.

\begin{lemma} \label{reconstrlem} 
Let $\D=(W,L,\{\Z b_\sigma\})$ be a $\Z$-root datum with associated
maximal torus normalizer $\nu(\D)$ and root subgroups $\nu(\D)_\sigma$. Then
$$
\{ x^2\, |\, x \in \nu(\D)_\sigma \setminus T_0^-(\sigma) \} =
\{h_\sigma\}
$$
for all reflections $\sigma\in W$. The analogous result holds for
$\Z_2$-root data.
\end{lemma}

\begin{proof}
If $x_1, x_2 \in \nu(\D)_\sigma \setminus
T_0^-(\sigma)$, then $x_2 = t x_1$ for some $t\in
T_0^-(\sigma)$ and hence $x_2^2 = (t x_1 t x_1^{-1}) x_1^2 =
(t \sigma(t)) x_1^2 = x_1^2$. On the other hand we have a commutative diagram
$$
\xymatrix{
1 \ar[r] & C_\sigma \oplus \Z\sigma \ar[r] \ar[d]^f &
\rho(W)_\sigma \ar[r] \ar[d] & \left<\sigma\right>
\ar[r] \ar@{=}[d] & 1\\
1 \ar[r] & T_0^-(\sigma) \ar[r] & \nu(\D)_\sigma \ar[r] &
\left<\sigma\right> \ar[r] & 1
}
$$
in which $Q_\sigma \subseteq \rho(W)_\sigma = \left<Q_\sigma\right>$
is mapped into $\nu(\D)_\sigma \setminus T_0^-(\sigma)$ by
Proposition~\ref{gen210}\eqref{proj}. Since $x^2 = \sigma$ for any
$x\in Q_\sigma$ and $f(\sigma) = h_\sigma$ this proves the
claim. Obviously the proof carries over verbatim to $\Z_2$-root data.
\end{proof}

We now enumerate the subgroups of $\nu(\D)$ which behave like the root
subgroups $\nu(\D)_\sigma$ in the sense of Lemma~\ref{reconstrlem}.

\begin{lemma} \label{torsorlemma}
Let $\D = (W,L,\{\Z b_\sigma\})$ be a $\Z$-root datum with corresponding
maximal torus normalizer $N = \nu(\D)$. Let $\sigma \in W$ be a
reflection, and let $N(\sigma)$ be the preimage under the projection
$N \to W$ of the subgroup $\left<\sigma\right>$. Denote by
$\mathcal{X}$  the set of subgroups $H$ of $N$ which sit in an exact sequence
\begin{equation} \label{Nsigmases}
1 \to T_0^-(\sigma) \to H \to \left< \sigma \right> \to 1
\end{equation}
and satisfy $\{x^2 \,|\, x \in H \setminus T_0^-(\sigma)\} =
\{h_\sigma\}$. Then $\mathcal{X}$ is a non-trivial
$H^1(\left<\sigma\right>;T)$-torsor with the action induced by the
obvious action of $\Der(\left<\sigma\right>,T) \subseteq
\Aut(N(\sigma))$ on $\mathcal{X}$.

The corresponding result for $\Z_2$-root data is also valid.
\end{lemma}

\begin{proof}
We already know that $N_\sigma\in \mathcal{X}$ by
Lemma~\ref{reconstrlem} so $\mathcal{X}$ is non-empty. 
Pick $x_0 \in N_\sigma \setminus T_0^-(\sigma)$; we have
$N_\sigma = \left< T_0^-(\sigma),x_0\right>$. It is clear that any
subgroup $H$ of $N$ which sits in an exact sequence of the form
\eqref{Nsigmases} must have the form $H=\left< T_0^-(\sigma),t
x_0\right>$ for some $t\in T$ and the computation 
$(tx_0)^2 = (t x_0 t x_0^{-1}) x_0^2 = t \sigma(t) h_\sigma$ shows
that $\left< T_0^-(\sigma),t x_0\right>\in \mathcal{X}$ if and only
if $t\in T^-(\sigma)$.

This shows that we have a transitive action of $T^-(\sigma)$ on
$\mathcal{X}$ given by $(t,\left< T_0^-(\sigma),x\right>) \mapsto
\left< T_0^-(\sigma),t x\right>$. Moreover it is easily seen that
under the identification $T^-(\sigma) = \Der(\left<\sigma\right>,T)$
this action corresponds to the natural action of
$\Der(\left<\sigma\right>,T) \subseteq \Aut(N(\sigma))$ on
$\mathcal{X}$. Since the stabilizer of a point in $\mathcal{X}$ equals
$T_0^-(\sigma) = \PDer(\left<\sigma\right>,T)$ we see that the
induced action of $H^1(\left<\sigma\right>;T)$ on $\mathcal{X}$ is
free and transitive, i.e.\ $\mathcal{X}$ is an
$H^1(\left<\sigma\right>;T)$-torsor.
\end{proof}

Next we have the following proposition, which will be used again in
Section~\ref{proof-sec}, in the proofs of Theorems~\ref{mainthmlie} and
\ref{mainthm2cg}.

\begin{prop} \label{H1calc}
Let $\D = (W,L,\{\Z b_\sigma\})$ be a $\Z$-root datum with corresponding
maximal torus normalizer $N = \nu(\D)$. Then $\Der(W,T) \cap
\Aut(N,\{N_\sigma\}) = \PDer(W,T)$. Similarly, if $\D =
(W,L,\{\Z_p b_\sigma\})$ is a $\Z_p$-root datum with corresponding
discrete maximal torus normalizer $\dN = \dnu(\D)$ then $\Der(W,\dT) \cap
\Aut(\dN,\{\dN_\sigma\}) = \PDer(W,\dT)$.
\end{prop}

The proof of this proposition relies on the following lemma
(compare \cite[Prop.~6.4]{matthey02}, \cite[Pf.~of~Lem.~2.5(i)]{HMS04}
and \cite[Pf.~of~Prop.~3.5]{JMO92}).

\begin{lemma} \label{Tsurjlemma}
Let $(W,L)$ be a finite $\Z$-reflection group and let $S$ be a system of
simple reflections. Then the homomorphism
$$
T \rightarrow \prod_{\sigma\in S} T_0^-(\sigma),\quad
t\mapsto \left(t \sigma(t)^{-1}\right)_{\sigma\in S}
$$
is surjective. The same result holds (replacing $T$ by $\dT$) when
$(W,L)$ is a finite $\Z_p$-reflection group of Coxeter type (i.e.,
$(W,L) \cong (W_1, L_1 \otimes_\Z \Z_p)$ for a $\Z$-reflection group
$(W_1, L_1)$).
\end{lemma}

\begin{proof}
Let $V=L\otimes_{\Z} \Q$. It suffices to prove that the homomorphism $\Psi
: V \to \prod_{\sigma\in S} V^-(\sigma)$ given by $\Psi(x) =
(x-\sigma(x))_{\sigma\in S}$ is surjective. The kernel of $\Psi$
obviously equals $V^W$ and since $\dim V - \dim V^W = |S|$ and $\dim
V^-(\sigma) = 1$ it follows that $\Psi$ is surjective.
\end{proof}

\begin{proof}[Proof of Proposition~\ref{H1calc}]
Assume first that $\D = (W,L,\{\Z b_\sigma\})$ is a $\Z$-root datum. Then
$f\in \Der(W,T)$ corresponds to the automorphism $\phi$ of $N$ given
by $\phi(x) = f(\bar{x}) x$. Since $\phi$ induces the identity on $W$
it follows that $\phi \in \Aut(N,\{N_\sigma\})$ if and only if
$\phi(N_\sigma) = N_\sigma$ for all $\sigma\in \Sigma$. By the properties of
$N_\sigma$ this is again equivalent to having $f(\sigma) \in
T_0^-(\sigma)$ for all $\sigma$. Hence we get
$\PDer(W,T) \subseteq \Der(W,T) \cap \Aut(N,\{N_\sigma\})$ since if
$f$ is a principal derivation given by $f(w) = (1-w)\cdot t$ for some
$t\in T$ we have $f(\sigma) = (1-\sigma)\cdot t \in (1-\sigma)T =
T_0^-(\sigma)$. To prove the reverse inclusion we have to prove
that any $f\in \Der(W,T)$ with $f(\sigma) \in T_0^-(\sigma)$ for all
$\sigma$ is a principal derivation. This follows directly from
Lemma~\ref{Tsurjlemma}.

Now let $\D = (W, L, \{\Z_p b_\sigma\})$ be a $\Z_p$-root datum. As above it
suffices to prove that any $f\in \Der(W,\dT)$ with $f(\sigma) \in
\dT_0^-(\sigma)$ for all $\sigma\in \Sigma$ is principal. For $p$ odd,
this follows from the fact that $H^1(W;\dT)=0$ by
\cite[Thm.~3.3]{kksa:thesis}. Now assume $p=2$. If $\D$ is of Coxeter type the
claim follows directly from Lemma~\ref{Tsurjlemma}. In the general
case we can write $\D$ as a product of a $\Z_2$-root datum of Coxeter
type and a number of copies of $\D_{\DI(4)}$, cf.\
\cite[Prop.~7.4]{DW05}. Since $\D_{\DI(4)}$ contains a subdatum isomorphic to
$\D_{\Spin(7)\twocom}$ (cf.\ \cite[Pf.~of~Prop.~9.12, $\DI(4)$
case]{DW05}) and $\left[W_{\DI(4)}:W_{\Spin(7)\twocom}\right] = 7$, it
follows that
$\D$ contains a subdatum $(W_1, L, \{\Z_2 b_\sigma\})$ of Coxeter type
such that $\left[W:W_1\right]$ is odd. From the above it now follows
that $[f]=0\in H^1(W_1;\dT)$. By a transfer argument the restriction
homomorphism $H^1(W;\dT) \to H^1(W_1;\dT)$ is injective proving the
claim in the general case.
\end{proof}

For a $\Z_p$-root datum $\D =(W,L,\{\Z_p b_\sigma \})$ we define the
{\em singular set}
\begin{equation} \label{singset}
S(\sigma) = \left< \dT^+_0(\sigma), h_\sigma \right>.
\end{equation}
One easily sees that this agrees with the definition in
\cite[Def.~7.3]{dw:center} (cf.\ the proof of
Proposition~\ref{reconstr} below), and that $S(\sigma) = \ker(\beta_\sigma \otimes_{\Z_p} \Z/p^\infty: \dT \to
\Z/p^\infty)$, which is analogous to the standard Lie theoretical definition.

For the proof of Theorem~\ref{sigmaident}, we need the following
result, which essentially summarizes elements of \cite[\S
9--10]{DW05}.

\begin{lemma} \label{centr-setup}
Let $X$ be a connected $2$-compact group with maximal discrete torus
$\dT$ and root datum $\D_X = (W, L, \{\Z_2 b_\sigma\}_{\sigma\in
\Sigma})$. Let $A$ be a subgroup of $\dT$, and let
$Y=\cC_X(A)_1$. Then $W_Y$ is a $\Z_2$-reflection group generated by
the reflections $\Sigma_Y =
\{\sigma\in \Sigma \,|\, A\subseteq S(\sigma)\}$ and $\D_Y$ identifies
with the subdatum $(W_Y, L, \{\Z_2 b_\sigma\}_{\sigma\in \Sigma_Y})$ of $\D_X$. 
Furthermore $\dN_Y$ is
the pull-back of $\dN_X$ and $\dnu(\D_Y)$ is the
pull-back of $\dnu(\D_X)$ along $W_Y \to W_X$.
\end{lemma}

\begin{proof}
The statement about the set of reflection $\Sigma_Y$ in $W_Y$ and the
result about normalizers follows from
\cite[Thm.~7.6(2)]{dw:center}. The fact that $\D_Y$ is a subdatum of
$\D_X$ now follows by definition of the root datum. The last statement
follows easily from Lemma~\ref{pullback}, cf. \cite[Lem.~10.1]{DW05}.
\end{proof}

Finally we need the following calculation.

\begin{lemma} \label{WDI4calc}
We have $H^1(W_{\DI(4)};\dT_{\DI(4)}) = \Z/2$ and
$H^1(\left<\sigma\right>;\dT_{\DI(4)}) = \Z/2$ for any reflection
$\sigma\in W_{\DI(4)}$.
\end{lemma}

\begin{proof}
Let $W=W_{\DI(4)}$ and $\dT=\dT_{\DI(4)}$ for short. Any reflection
$\sigma\in W$ acts non-trivially on ${}_2 \dT$, so up to conjugation the action of
$\sigma$ is given by $(t_1,t_2,t_3)\mapsto (t_2,t_1,t_3)$. It now
follows directly that $H^1(\left<\sigma\right>;\dT) \cong
\dT^-(\sigma)/\dT_0^-(\sigma) = \Z/2$.

To see the first claim note that $W\cong Z(W) \times
\GL_3(\F_2)$, where $Z(W) = \left<-1\right>$. The associated Lyndon-Hochschild-Serre
spectral sequence for computing $H^*(W;\dT)$ has $E_2$-term
$$
E_2^{s,t} = H^s(\GL_3(\F_2);H^t(\left<-1\right>;\dT)).
$$
For $t$ odd we have $E_2^{s,t}=0$ since
$H^t(\left<-1\right>;\dT)=0$. For position reasons we thus get
$H^1(W;\dT) \cong E_2^{1,0}$. Since $H^0(\left<-1\right>;\dT) \cong
(\F_2)^3$ with the natural action of $\GL_3(\F_2)$, we get $H^1(W;\dT) \cong
H^1(\GL_3(\F_2);{\F_2}^3)$. It is well known that
$H^1(\GL_3(\F_2);{\F_2}^3) = \Z/2$; probably the easiest way to see
this directly is to observe that for $P_{\text{triv}}$, the projective cover
of the trivial $\F_2[\GL_3(\F_2)]$-module, the second radical (Loewy) layer
$\rad(P_{\text{triv}})/\rad^2(P_{\text{triv}})$ equals the standard representation
${\F_2}^3$ plus its dual by \cite[p.~216]{benson84} (or a direct
calculation), so the result follows from the minimal resolution; see
also \cite[Prop.~4(a)]{sah74}, \cite[Table~C]{JP76} or \cite[Table~I]{bell78a}. 
\end{proof}

\begin{proof}[Proof of Theorem~\ref{sigmaident}]
In the Lie group case \eqref{maplie}, the existence of the isomorphism
follows from \cite[Thm.~4.4]{tits66} after translating the definitions
(isomorphism in \cite[Thm.~4.4]{tits66} means ``isomorphism in the category
$\mathscr N$'' defined in \cite[3.1]{tits66} which exactly means
sending $\nu(\D)_\sigma$ to $N_G(T)_\sigma$; Tits states his results
for reductive algebraic groups, which translate into results for
compact Lie groups via the usual method). A different exposition of this 
is given in \cite[Pf.~of~Thm.~5.7]{DW05}
noting that, in the notation of \cite{DW05}, $q_i$ is sent to $x_i$
and hence $\nu(\D)_{\sigma_i} = \left< T_0^-(\sigma_i),q_i\right>$ is
sent to $N_G(T)_{\sigma_i} = \left<  T_0^-(\sigma_i),x_i\right>$ for
all simple reflections $\sigma_i$, so $\nu(\D)_\sigma$ is sent
to $N_G(T)_\sigma$ for all $\sigma$ since both the topological and the
algebraic subgroups are well behaved under conjugation. The
uniqueness follows directly from Proposition~\ref{H1calc}.

We now prove part \eqref{mappcg}. Assume first that $p>2$. In this
case the result follows immediately from the definitions above: There
is an isomorphism of extensions $\dN_X \cong \dT\rtimes W =
\dnu(\D_X)$ by \cite[Thm.~1.2]{kksa:thesis} and by definition the
topological root subgroups $(\dN_X)_\sigma$ are the images of the
algebraic root subgroups $\dnu(\D_X)_\sigma = \dT^-(\sigma) \rtimes
\left<\sigma\right>$ under this isomorphism. Since $H^1(W;\dT) = 0$
\cite[Thm.~3.3]{kksa:thesis}, uniqueness also follows. 

Now assume that $p=2$. Since the uniqueness part follows directly from
Proposition~\ref{H1calc}, we only need to prove the existence of the
isomorphism. First note that if $X$ is the $2$-completion of a compact
connected Lie group then  the result follows from part~\eqref{maplie}
since in this case the definitions and constructions for $2$-compact
groups agree in the obvious way with those for compact Lie groups
(cf.\ \cite[Lems.~8.1 and 9.17]{DW05}). In particular the result
holds when $X$ has semi-simple rank $1$, by the classification
mentioned in the beginning of this section. The general case can be
reduced to the case where $X$ has semi-simple rank $2$. We therefore
first explain this case (where we use some group cohomology
calculations, but essentially avoid classification results for $2$-compact groups)
and then explain how the general case
follows from this.

Assume that $X$ has semi-simple rank $2$. Since $\D_{\DI(4)}$ has
semi-simple rank $3$, \cite[Prop.~7.4]{DW05} shows that $\D_X \cong
\D_G \otimes_\Z \Z_2$ for some compact connected Lie group $G$ of
semi-simple rank $2$. In particular $G$ has Cartan type $A_1A_1$,
$A_2$, $B_2 = C_2$ or $G_2$.

If $G$ has Cartan type $A_1A_1$, then $X' = X/\cZ(X)$ is a connected
centerfree $2$-compact group of rank $2$ and $(W_{X'},L_{X'} \otimes \Q)$
is reducible. Hence by \cite[Thm.~1.3]{dw:split} $X'$ splits as
a product of two center-free $2$-compact groups of rank $1$. Hence by
the classification of $2$-compact groups of rank $1$, $X' \cong
\SO(3)\twocom \times \SO(3)\twocom$, and consequently $X$ is the
$2$-completion of a compact connected Lie group, and the result
follows.

In the case where $G$ has Cartan type $A_2$, $B_2$ or $G_2$, an inspection
of the classification of $\Z$-root data reveals that $\D_X
\cong \D_G\otimes_{\Z} \Z_2$ is isomorphic to $\D_{H \times (S^1)^n}
\otimes_{\Z} \Z_2$, $n\geq 0$, for $H=\SU(3)$, $\Spin(5)$,
$\SO(5)$, $(\Spin(5)\times S^1)/C_2$ (with $C_2$ the
unique ``diagonally'' embedded central subgroup of order $2$ in
$\Spin(5)\times S^1$), or $G_2$. Since $\D_X \cong \D_{H\twocom} \times
\D_{((S^1)^n)\twocom}$, \cite[Thm.~1.4]{dw:split} implies that $X
\cong X' \times ((S^1)^n)\twocom$ for a connected $2$-compact group
$X'$ with $\D_{X'} \cong \D_H \otimes_{\Z} \Z_2$. It is easily seen from the
definitions that part~\eqref{mappcg} will hold for $X$ if it holds for
$X'$. We now prove this by going through the five possibilities for $H$
one by one:

For $H=\SU(3)$ and $G_2$ we have $H^2(W_{X'};\dT_{X'}) = 0$, so
the extensions $\dnu(\D_{X'})$ and $\dN_{X'}$ are
isomorphic. Moreover, as $H^1(\left<\sigma\right>;\dT_{X'})=0$ for any
reflection $\sigma\in W_{X'}$, Lemma~\ref{torsorlemma} shows that any
isomorphism between the two extensions will preserve root subgroups.

In the case $H=\Spin(5)$, there is, up to conjugation, a unique non-central
element of order two $\nu: B\Z/2 \to BX'$ in $X'$. By the description
of the root datum of a centralizer, Lemma~\ref{centr-setup},
$Y=\cC_{X'}(\nu)$ is a connected $2$-compact group with $\D_Y \cong \D_{(\SU(2)\times
\SU(2))\twocom}$. In particular the extensions $\dnu(\D_Y)$ and
$\dN_Y$ are isomorphic by the result in the $A_1A_1$ case, and
moreover a computation shows that the restriction $H^2(W_{X'};\dT) \cong (\Z/2)^2 \to
H^2(W_Y;\dT) \cong (\Z/2)^4$ is injective. Combining this with
Lemma~\ref{centr-setup} shows that the extensions $\dnu(\D_{X'})$ and
$\dN_{X'}$ are isomorphic. We now argue that this isomorphism can be
taken to preserve the root subgroups:
Let $f\in \Der(W_{X'},\dT)$. By the proof of Proposition~\ref{H1calc} it
follows that if $f|_{\left<\sigma\right>}\in
\PDer(\left<\sigma\right>,\dT)$ for all reflections $\sigma$, then
$f\in \PDer(W_{X'},\dT)$. In other words the intersection of the
kernels of the restrictions, $\bigcap_\sigma \ker\left(H^1(W_{X'};\dT)
\to H^1(\left<\sigma\right>;\dT)\right)$, equals $0$, where the intersection
runs over all reflections $\sigma$. However the kernel of the
restriction only depends on the conjugacy class of $\sigma$ (cf.\
e.g., \cite[Ex.~III.9.1]{brown94}), so letting $\sigma_1$ and
$\sigma_2$ denote representatives for the two conjugacy classes of
reflections in $W_{X'}$, we have
$$
\ker\left(H^1(W_{X'};\dT) \to H^1(\left<\sigma_1\right>;\dT)\right)
\cap \ker\left(H^1(W_{X'};\dT) \to
H^1(\left<\sigma_2\right>;\dT)\right) = 0.
$$
A direct computation shows that (up to exchange of $\sigma_1$ and
$\sigma_2$), we have $H^1(W_{X'};\dT) = \Z/2$,
$H^1(\left<\sigma_1\right>;\dT) = \Z/2$ and
$H^1(\left<\sigma_2\right>;\dT) = 0$, so we conclude that the
restriction $H^1(W_{X'};\dT) \to H^1(\left<\sigma_1\right>;\dT)$ is
surjective. Hence, by Lemma~\ref{torsorlemma}, $H^1(W;\dT)$ acts
transitively on the subgroups of $\dnu(\D_{X'})$ which look
like $\dnu(\D_{X'})_{\sigma_1}$. Therefore any isomorphism between
$\dnu(\D_{X'})$ and $\dN_{X'}$ can be modified by an element of
$\Der(W_{X'},\dT) \subseteq \Aut(\dnu(\D_{X'}))$ in such a way that
$\dnu(\D_{X'})_{\sigma_1}$ is sent to
$(\dN_{X'})_{\sigma_1}$. Moreover as $H^1(\left<\sigma_2\right>;\dT) =
0$, Lemma~\ref{torsorlemma} shows that this modified isomorphism
automatically will send $\dnu(\D_{X'})_{\sigma_2}$ to
$(\dN_{X'})_{\sigma_2}$. Since any reflection is conjugate to
$\sigma_1$ or $\sigma_2$ and the root subgroups are well behaved
under conjugation, the modified isomorphism will send every root
subgroup of $\dnu(\D_{X'})$ to the corresponding root subgroup of
$\dN_{X'}$ as desired.

When $H=\SO(5)$, $W_{X'}$ has two conjugacy classes of reflections
represented by $\sigma_1$ and $\sigma_2$ say. Moreover $X'$ has two
conjugacy classes of elements of order two, $\nu_1, \nu_2:
B\Z/2 \to BX'$, and we find that (up to permutation)
$Y_1 = \cC_{X'}(\nu_1)_1$ has root datum isomorphic to
$\D_{\SO(4)\twocom}$ and that $Y_2 = \cC_{X'}(\nu_2)_1$ satisfies
$\D_{Y_2} \cong \D_{(\SO(3)\times S^1)\twocom}$. Hence the extensions
$\dnu(\D_{Y_i})$ and $\dN_{Y_i}$ are isomorphic for $i=1$ and $2$ by
the result in the $A_1A_1$ case and the semi-simple rank one case.
Furthermore a calculation shows that the homomorphism
$$
H^2(W_{X'};\dT) \cong (\Z/2)^2 \to H^2(W_{Y_1};\dT) \times H^2(W_{Y_2};\dT) \cong \Z/2
\times \Z/2
$$
given by restrictions is an isomorphism, so by Lemma~\ref{centr-setup} the
extensions $\dnu(\D_{X'})$ and $\dN_{X'}$ are isomorphic. Since the
Weyl group action is the same as in the case $H=\Spin(5)$, the
argument given in that case shows that we can modify the isomorphism
between $\dnu(\D_{X'})$ and $\dN_{X'}$ by an element in
$\Der(W_{X'},\dT) \subseteq \Aut(\dnu(\D_{X'}))$ to obtain an
isomorphism which preserves root subgroups.

Finally, when $H=(\Spin(5)\times S^1)/C_2$, there is an element of
order two $\nu:B\Z/2 \to BX'$ such that $Y=\cC_{X'}(\nu)$ is a
connected $2$-compact group with $\D_Y \cong \D_{H'\twocom}$, where
$H' = (\SU(2)\times \SU(2) \times S^1)/C_2$, where $C_2$ is the unique
``diagonally'' embedded central subgroup of order $2$ in $\SU(2)\times
\SU(2)\times S^1$. In particular $\dnu(\D_Y)$ and $\dN_Y$ are
isomorphic by the above and since a computation shows that the
restriction $H^2(W_{X'};\dT) \cong \Z/2 \to H^2(W_Y;\dT) \cong
(\Z/2)^2$ is injective we conclude that $\dnu(\D_{X'})$ and
$\dN_{X'}$ are isomorphic. Since $W_{X'}$ has two conjugacy classes of
reflections represented by $\sigma_1$ and $\sigma_2$, the computations
$H^1(W_{X'};\dT) = (\Z/2)^2$, $H^1(\left<\sigma_1\right>;\dT) = \Z/2$ and
$H^1(\left<\sigma_2\right>;\dT) = \Z/2$ combined with the argument in the
case $H=\Spin(5)$ implies that the restriction $H^1(W_{X'};\dT) \to
H^1(\left<\sigma_1\right>;\dT) \times H^1(\left<\sigma_2\right>;\dT)$
is an isomorphism, so as above we can choose an isomorphism between
$\dnu(\D_{X'})$ and $\dN_{X'}$ which preserves root subgroups.
This finishes the proof of part~\eqref{mappcg} in the case where $X$
has semi-simple rank at most $2$. 

We now proceed to give the proof in the general
case. We follow \cite[Pf.~of~Prop.~9.12]{DW05} closely, but since extra arguments
are needed we choose to continue in some detail.
By \cite[Thm.~6.1]{dw:split} and \cite[Prop.~7.4]{DW05} we have $X
\cong X_1 \times X_2$ where $\D_{X_1}$ is of Coxeter type and
$\D_{X_2}$ is a product of a number of copies of $\D_{\DI(4)}$, and it
is hence enough to treat these two cases separately.

If $\D_X$ is of Coxeter type, we let $\{\sigma_i\}_{i\in I}$ be a set of
simple reflections. Then by \cite[\S 2]{tits66} (cf. also
\cite[4.6 and Def.~6.15]{DW05}) $\dnu(\D_X)$ is generated by $\dT$ and symbols
$q_i$, $i\in I$, subject to the relations $q_i^2 = h_{\sigma_i}$, $q_i
t q_i^{-1} = \sigma_i(t)$ for $t\in \dT$ and
\begin{equation}
\underbrace{q_i q_j \ldots}_{\text{$m_{ij}$ terms}} =
\underbrace{q_j q_i \ldots}_{\text{$m_{ij}$ terms}}
\end{equation}
for $i\neq j$, where $m_{ij}$ denotes the order of $\sigma_i
\sigma_j$. Choose elements $x_i \in (\dN_X)_{\sigma_i}\setminus
\dT^-_0(\sigma_i)$ for each $i\in I$. As in \cite[Pf.~of~Prop.~9.12,
Coxeter case]{DW05} we see that the elements $x_i$ satisfies the first
two relations and that it suffices to check the third relation for the
subgroup $Y = \cC_X(\dT_0^{\langle \sigma_i,\sigma_j\rangle})$ of $X$,
which has semi-simple rank $2$. In this case there is a compact
connected Lie group $G$ with $\D_Y \cong \D_{G\twocom}$. By the above,
part~\eqref{mappcg} of the theorem holds for $Y$, so there is an
isomorphism between the extensions $\dN_Y$ and $\dN_{G\twocom}$ which
preserves root subgroups. Hence it suffices to check the third
relation for $G\twocom$ where it follows from
the classical result of Tits (cf.\ \cite[Thm.~5.8 and Lem.~9.17]{DW05}).
Hence there is an isomorphism $\dnu(\D_X)
\xrightarrow{\cong} \dN_X$ induced by $q_i \mapsto x_i$ which clearly
preserves root subgroups. This proves the result in the Coxeter case.

Finally assume that $\D_X \cong \D_{\DI(4)}$. In this case there is a unique
element of order two $\nu: B\Z/2 \to BX$ and
$Y=\cC_X(\nu)$ satisfies $\D_Y \cong \D_{\Spin(7)\twocom}$. Hence there is
an isomorphism between the extensions $\dnu(\D_Y)$ and $\dN_Y$ by the
result in the Coxeter case. As $\left[W_X:W_Y\right] = 7$ is odd, the
restriction $H^2(W_X;\dT) \to H^2(W_Y;\dT)$  is injective, so
Lemma~\ref{centr-setup} shows that $\dnu(\D_X)$ and $\dN_X$ are
isomorphic. Let $\sigma\in W = W_{\DI(4)}$ be a reflection. Since all
reflections in $W$ are conjugate, the argument given above shows that
the restriction $H^1(W;\dT) \to H^1(\left<\sigma\right>;\dT)$ is
injective. Now $H^1(W;\dT) \cong \Z/2$ and
$H^1(\left<\sigma\right>;\dT) \cong \Z/2$ by Lemma~\ref{WDI4calc}, so
the restriction is an isomorphism. As above this implies that we can
modify the isomorphism between $\dnu(\D_X)$ and $\dN_X$ to an
isomorphism which sends root subgroups to root subgroups as desired.
\end{proof}


\section{Proofs of Theorems~\ref{mainthmlie} and \ref{mainthm2cg}} \label{proof-sec}

In this section we prove Theorems~\ref{mainthmlie}
and \ref{mainthm2cg} using the material from the previous
sections. First we prove the following result which says that if
$\phi$ is an automorphism of $W$ which sends reflections to
reflections, then the class in $H^2(W;\Z[\Sigma])$ corresponding to
the reflection extension $\rho(W)$ is fixed under $(\phi,\phi^{-1})^*:
H^2(W;\Z[\Sigma]) \to H^2(W;\Z[\Sigma])$, where $\phi^{-1}$ induces a
map $\Z[\Sigma] \to \Z[\Sigma]$ by assumption.

\begin{lemma} \label{rhoWauto}
Let $W$ be a finite $\Q_2$-reflection group with associated reflection
extension $1 \to \Z[\Sigma] \xrightarrow{\iota} \rho(W) \xrightarrow{\pi} W \to
1$. For any $\phi\in \Aut(W)$ with $\phi(\Sigma)=\Sigma$, there is an
automorphism $\psi$ of $\rho(W)$ fitting into the commutative diagram
$$
\xymatrix@C=40pt{
1 \ar[r] & \Z[\Sigma] \ar[r]^-{\iota} \ar@{=}[d] & \rho(W) \ar[r]^-{\pi}
\ar[d]^{\psi} & W \ar[r] \ar@{=}[d] & 1 \\
1 \ar[r] & \Z[\Sigma] \ar[r]^-{\iota\circ\phi} & \rho(W)
\ar[r]^-{\phi^{-1}\circ\pi} & W \ar[r] & 1.
}
$$
Moreover $\psi$ is unique up to conjugation by an element in $\Z[\Sigma]$.
\end{lemma}

\begin{proof}
Define a category $\mathscr{D}$ as follows
(cf.\ \cite[Ch.~III, \S 8]{brown94}): The objects are pairs $(G,M)$, where
$G$ is a group and $M$ is a $G$-module, and a morphism from $(G,M)$ to
$(G',M')$ is a pair of maps $(\alpha:G \to G', f:M' \to M)$ where $\alpha$ is a group homomorphism and $f$ is a
$G$-module homomorphism, i.e.\ $f(\alpha(g)\cdot m') = g\cdot
f(m')$ for $g\in G$ and $m'\in M'$. On the level of cohomology this
induces a homomorphism $(\alpha,f)^* : H^*(G';M') \to H^*(G;M)$, which
in degree $2$ takes the equivalence class of an extension $0 \to M'
\to E' \to G' \to 1$ to the equivalence class of the extension $0 \to M
\to E \to G \to 1$ obtained by taking the pullback along $\alpha:G
\to G'$ followed by the pushforward along $f:M' \to M$.

Now let $k\in H^2(W;\Z[\Sigma])$ denote the class of the reflection
extension. The class of the extension $0 \to \Z[\Sigma]
\xrightarrow{\iota\circ\phi} \rho(W) \xrightarrow{\phi^{-1}\circ\pi} W
\to 1$ is then given by $(\phi,\phi^{-1})^*(k)$ and the first part of the
proposition is that this equals $k$, since in that case $(\phi,\phi^{-1})$
induces an automorphism $\psi: \rho(W) \to \rho(W)$ as wanted.

Let $\Sigma = \bigcup_i \Sigma_i$ denote the partition of the
set of reflections $\Sigma$ in $W$ into conjugacy classes. By
construction $k\in H^2(W;\Z[\Sigma])$ is given by a sum of elements
$k_i\in H^2(W;\Z[\Sigma_i])$. The map $(\phi,\phi^{-1})$ also induces
homomorphisms $H^2(W;\Z[\Sigma_i]) \to H^2(W;\Z[\phi^{-1}(\Sigma_i)])$
and we prove below that
$(\phi,\phi^{-1})^*(k_i) = k_{i'}$ where $\Sigma_{i'} =
\phi^{-1}(\Sigma_i)$. Thus $(\phi,\phi^{-1})^*$ permutes the $k_i$'s
and hence fixes $k$.

To see this claim, let $\tau_i\in \Sigma_i$ and define $C_i$ and
$C_i^\perp$ as in Section~\ref{reflext-sec}. Moreover, let
$\text{incl}_i:C_i \to W$ denote the inclusion, $p_i: C_i =
\left<\tau_i\right> \times C_i^\perp \to \left<\tau_i\right>$ the
projection and $f_i: \Z[W/C_i] \to \Z$ the $C_i$-module homomorphism
given by
$$
f_i(wC_i) =
\begin{cases}
1 & \text{if $w\in C_i$}, \\
0& \text{if $w\notin C_i$}.
\end{cases}
$$
By \cite[Ex.~III.8.2]{brown94} (see also \cite[Lem.~4.8]{DW05})
the isomorphism $H^2(W;\Z[W/C_i]) \xrightarrow{\cong} H^2(C_i;\Z)$ from
Shapiro's lemma is induced by $(\text{incl}_i,f_i)$. Choosing $\tau_{i'} =
\phi^{-1}(\tau_i)$ as the representative from $\Sigma_{i'} =
\phi^{-1}(\Sigma_i)$, the claim that $(\phi,\phi^{-1})^*(k_i) =
k_{i'}$ then follows from the commutativity of the diagram
$$
\xymatrix@C=42pt{
H^2(\left<\tau_i\right>;\Z) \ar[r]^-{p_i^*} \ar[d]_-{\phi^*} & H^2(C_i;\Z) &
H^2(W;\Z[W/C_i]) \ar[l]^-{\cong}_-{(\text{incl}_i,f_i)^*} \ar@{=}[r] &
H^2(W;\Z[\Sigma_i]) \ar[d]^-{(\phi,\phi^{-1})^*} \\
H^2(\left<\tau_{i'}\right>;\Z) \ar[r]^-{p_{i'}^*} & H^2(C_{i'};\Z) &
H^2(W;\Z[W/C_{i'}]) \ar[l]^-{\cong}_-{(\text{incl}_{i'},f_{i'})^*} \ar@{=}[r] &
H^2(W;\Z[\Sigma_{i'}])
}
$$
which is easily established by checking that the corresponding diagram
in $\mathscr{D}$ commutes. To see the second part, note first that
$$
H^1(W;\Z[\Sigma]) = \bigoplus_i H^1(W;\Z[\Sigma_i]) \cong \bigoplus_i
H^1(C_i;\Z) = 0.
$$
It is clear that $\psi$ is unique up to an automorphism of $\rho(W)$
inducing the identity on $\Z[\Sigma]$ and on $W$. If $\chi$ is any
such automorphism, we get a well-defined derivation $f:W \to
\Z[\Sigma]$ given by $f(w) = \chi(\widetilde{w}) {\widetilde{w}}^{-1}$, where
$\widetilde{w}\in \rho(W)$ is any lift of $w\in W$. Since
$H^1(W;\Z[\Sigma])=0$, this must be a principal derivation which means
that $\chi$ must be conjugation by an element in $\Z[\Sigma]$, so
$\psi$ is unique up to conjugation by an element in $\Z[\Sigma]$.
\end{proof}

The next proposition, which we will need later on,
shows how to recover the root datum $\D$ from the
isomorphism type of the associated maximal torus normalizer
$\nu(\D)$, detailing a remark of Dwyer-Wilkerson \cite[Rem.~5.6]{DW05}.
A result of this type was first proved via different arguments by
Curtis-Wiederhold-Williams \cite{CWW74}, for $G$ a compact connected
Lie group with finite center, and by
Osse \cite{osse97} and Notbohm \cite{notbohm95} for $G$ an arbitrary
compact connected Lie group.

First note that $T$ is characterized as the largest divisible subgroup
in $\nu(\D)$, and therefore the extension $1 \to T \to \nu(\D) \to W
\to 1$ can be recovered, so the only issue is how to recover the
markings $h_\sigma$.

\begin{prop}[{\cite[Rem.~5.6]{DW05}}] \label{reconstr}
Let $\D = (W,L,\{\Z b_\sigma\})$ be a $\Z$-root datum and let $N = \nu(\D)$ be
the associated maximal torus normalizer. Let $\sigma\in W$ be a reflection and
let $N(\sigma)$ denote the pullback of $1 \to T \to N \to W \to 1$
along the inclusion $\left<\sigma\right> \to W$. Then the element
$h_\sigma$ is the unique element in $T_0^-(\sigma) \cap \{ x^2\, |\, x\in
N(\sigma) \setminus T \}$ which is a marking for the reflection
$\sigma$ on $T$ (i.e.\ an element $x \in T_0^-(\sigma)$ such that
$x^2=1$ and $x \neq 1$ if $\sigma$ acts non-trivially on ${}_2 T$,
cf.\ \cite[Def.~2.12]{DW05}).

For $\Z_2$-root data the corresponding result (obtained by replacing
$T$ by $\dT$ and $\nu(\D)$ by $\dnu(\D)$) also holds.
\end{prop}

\begin{proof}
By Lemma~\ref{reconstrlem}, $h_\sigma$ lies in $T_0^-(\sigma) \cap \{ x^2\, |\, x\in
N(\sigma) \setminus T \}$. We need to see that there are no other
elements in this set which are markings for $\sigma$.
This is clear when $\sigma$ acts non-trivially on ${}_2T$,
since $h_\sigma$ is then by definition the unique marking for $\sigma$
in $T_0^-(\sigma)$. In general let $x\in
N(\sigma)\setminus T$ with $x^2=h_\sigma$. For $t\in T$ we have $(t x)^2 = (t x
t x^{-1}) x^2 = t \sigma(t) h_\sigma$ proving that $\{ x^2\, |\, x\in
N(\sigma) \setminus T \} = h_\sigma \cdot (1+\sigma)(T) = h_\sigma
T_0^+(\sigma)$. Hence our claim is that $h_\sigma$ is the only marking
for $\sigma$ in $T_0^-(\sigma) \cap h_\sigma T_0^+(\sigma) = h_\sigma
(T_0^-(\sigma) \cap T_0^+(\sigma))$. If $\sigma$ acts trivially on
${}_2 T$, the action of $\sigma$ is (up to an automorphism of $T$) given by $(t_1, t_2,
\ldots, t_n) \mapsto (t_1^{-1}, t_2, \ldots, t_n)$ so
$T_0^-(\sigma) \cap T_0^+(\sigma) = 1$ proving the claim in this
case also.

This proves the first part of the proposition. The argument for
$\Z_2$-root data is identical.
\end{proof}

We are now ready to prove Theorem~\ref{mainthmlie}.

\begin{proof}[Proof of Theorem~\ref{mainthmlie}]
We divide the proof of the first part into several
steps.

{\em Step $1$:} The restriction homomorphism $\Aut(N) \to N_{\Aut(T)}(W)$ has
image contained in $\Aut(\D)$:
We just need to see that for any $\phi \in \Aut(N)$ we have
$\phi(h_\sigma) = h_{\phi\sigma\phi^{-1}}$. However this follows
from the description of $h_\sigma$ in terms of $N$ as given in
Proposition~\ref{reconstr} above.

{\em Step $2$:} The restriction $\Aut(N) \to \Aut(\D)$ is surjective:
Clearly any automorphism $\phi \in
\Aut(\D)$ induces an automorphism of $W$ preserving $\Sigma$, which we
(by abuse of notation) also denote $\phi$. By
Lemma~\ref{rhoWauto} there thus exists an automorphism $\psi$ of
$\rho(W)$ satisfying $\pi\circ\psi = \phi\circ\pi$ and $\psi\circ\iota
= \iota\circ\phi$. The first condition shows that
$(\phi,\psi)$ defines an automorphism of $T\rtimes \rho(W)$ and the
second that the subgroup $\left< (h_\sigma,\sigma^{-1})\, |\,
\sigma\in\Sigma \right>$ is preserved by this automorphism (here
$\sigma^{-1}$ denotes the inverse of $\sigma\in \Z[\Sigma]\subseteq \rho(W)$). Hence
$(\phi,\psi)$ gives a well-defined automorphism of $N = (T\rtimes
\rho(W))/ \left< (h_\sigma,\sigma^{-1})\, |\, \sigma\in\Sigma \right>$
which obviously restricts to $\phi$ on $T$. 
Picking a $\psi$ for each $\phi \in \Aut(\D)$ gives
an explicit set
theoretical splitting $s:\Aut(\D) \to \Aut(N)$. (Note that $s$ depends
on the choices of $\psi$.)

{\em Step $3$:} The kernel of the restriction homomorphism $\Aut(N) \to
\Aut(\D)$ equals the group of derivations $\Der(W,T)$: First note that
we may view $\Der(W,T)$ as a subgroup of $\Aut(N)$ by sending $f \in
\Der(W,T)$ to $\phi_f \in \Aut(N)$ given by $\phi_f(x) = f(\bar
x)x$. Likewise if $\phi \in \Aut(N)$ restricts to the identity on $T$
then we may define $f \in \Der(W,T)$ by $f(w) = \phi(\widetilde{w})
\widetilde{w}^{-1}$, where $\widetilde{w} \in N$ is a lift of $w \in
W$, and with this definition $\phi_f = \phi$.

The three claims above show that we have a short exact sequence
\begin{equation} \label{autNexactseq}
1 \to \Der(W,T) \to \Aut(N) \to \Aut(\D) \to 1.
\end{equation}
Now consider the subgroup $\Inn(N)$ of $\Aut(N)$. The
image under the restriction homomorphism $\Inn(N) \to N_{\Aut(T)}(W)$
obviously equals $W$. Moreover it is easily seen that under the
identification from Step $3$, the kernel equals the group $\PDer(W,T)$
of principal derivations. Hence \eqref{autNexactseq} has the exact
subsequence $1 \to \PDer(W,T) \to \Inn(N) \to W \to 1$ and
the quotient exact sequence is the exact sequence in the
theorem.

{\em Step $4$:} The composition $\Aut(\D) \xrightarrow{s} \Aut(N) \to \Out(N)$
is a well-defined homomorphism: From the construction of $s$
and Lemma~\ref{rhoWauto} it follows that the value of $s(\phi)$
is well-defined up to an automorphism which on $T\rtimes \rho(W)$ has
the form $(\id_T,c_a)$ where $c_a$ denotes conjugation in $\rho(W)$ by
an element $a\in \Z[\Sigma]$. Obviously this agrees with conjugation
by the element $(1,a) \in T\rtimes \rho(W)$. Thus the composition
$\Aut(\D) \xrightarrow{s} \Aut(N) \to \Out(N)$ is independent of the choices
involved and by construction it is clear that the composition is a
homomorphism.

{\em Step $5$:} The homomorphism $\Aut(\D) \to \Out(N)$ from Step $4$
factors through $\Out(\D)$: We have to see that $s$ sends elements of
$W$ to elements of $\Inn(N)$. For $w\in W$ choose $x\in \rho(W)$ with
$\pi(x) = w$. It is easily checked that we can take $\psi = c_x$
(i.e., conjugation by $x$) in the definition of $s(w)$. Then $s(w)$ is
induced by the automorphism of $T\rtimes \rho(W)$ given
by $(t\mapsto w\cdot t,c_x)$ and since this
automorphism agrees with conjugation by $(1,x)\in T\rtimes \rho(W)$ the
claim follows.

{\em Step $6$:} The values of $s$ belong to $\Aut(N,\{N_\sigma\})$:
Let $\phi\in\Aut(\D)$ and choose $\psi\in \Aut(\rho(W))$ as in Step
$2$. Since $\psi$ stabilizes $\Z[\Sigma]$ as a set and agrees with
$\phi$ on this subgroup, it follows that $\psi(Q_\sigma) =
Q_{\phi\sigma\phi^{-1}}$. Obviously $\phi(T_0^-(\sigma)) =
T_0^-(\phi\sigma\phi^{-1})$ and hence $s(N_\sigma) =
N_{\phi\sigma\phi^{-1}}$ by definition.

{\em Step $7$:} By now we know that the splitting $s:\Out(\D) \to \Out(N)$
takes values in $\Out(N,\{N_\sigma\})$. It then follows from Proposition~\ref{H1calc}
 that the image of $s$ equals $\Out(N,\{N_\sigma\})$ as
claimed. This finishes the proof of the first part of the theorem.

Assume finally that $G$ is a compact connected Lie group with root
datum $\D_G$. By Theorem~\ref{sigmaident}\eqref{maplie} we can fix an identification
of $N_G(T)$ with $N = \nu(\D_G)$ in such a way that the subgroups
$N_G(T)_\sigma$ and $N_\sigma$ correspond to each other. 

Let $\Aut(G,T)$ denote the subgroup of $\Aut(G)$ which sends $T$ to $T$. The
composition $\Aut(G,T) \xrightarrow{\res} \Aut(N) \to \Out(N)$ obviously
factors through $\Out(G)$. Moreover by definition of the
subgroups $N_G(T)_\sigma$
(cf.\ Section~\ref{normalext-sec} and \cite[Def.~5.1]{DW05}) it is clear
that this composition takes values in $\Out(N,\{N_\sigma\})$.
The homomorphism $\Aut(G,T) \to \Out(G)$ is surjective since all
maximal tori in $G$ are conjugate. Thus we have a well-defined homomorphism
$\Out(G) \to \Out(N,\{N_\sigma\})$. Similarly using the
definition of the elements $h_\sigma$  (cf.\
Section~\ref{normalext-sec} and \cite[Def.~5.3]{DW05}), we get a
well-defined restriction homomorphism $\Out(G) \to \Out(\D_G)$. It is
now clear that the diagram
$$
\xymatrix{
 & \Out(N,\{N_\sigma\}) \ar[dr]^{\res} & \\
\Out(G) \ar[ur]^{\res} \ar[rr]^{\res} & & \Out(\D_G)
}
$$
commutes. The homomorphism $\Out(N,\{N_\sigma\}) \to \Out(\D_G)$ is an
isomorphism by the above and it is a fundamental theorem in Lie theory that
$\Out(G) \to \Out(\D_G)$ is an isomorphism
(cf.\ e.g.\ \cite[\S 4, no.~10, Prop.~18]{bo9}). This proves the theorem.
\end{proof}

\begin{rem}
Note that e.g., for $G = F_4$ the sequence $1 \to \Der(W,T) \to \Aut(N)
\to \Aut(\D) \to 1$ is {\em not} split:  By \cite{HMS04} $H^1(W_{F_4};T) =
0$ and clearly $\Out(\D) = 1$ (cf.\ \cite[\S 12.2, Table
1]{humphreys78}). Since $T^W = 0$, the sequence identifies with $1 \to T
\to N_{F_4}(T) \to W_{F_4} \to 1$, which is not split by \cite[Section
4, Application I]{CWW74}.
\end{rem}

\begin{lemma} \label{autlemmapodd}
Let $\D$ be a $\Z_p$-root datum, $p$ odd, and let $\dN = \dnu(\D)$. Then
$\Aut(\dN,\{\dN_\sigma\}) = \Aut(\dN)$.
\end{lemma}

\begin{proof}
Recall that $\dN = \dT \rtimes W$ by definition.
Let $\phi \in \Aut(\dN)$. Since $\dT$ is a characteristic subgroup
of $\dN$, $\phi$ induces an automorphism $\overline{\phi} : W
\to W$. Defining $\phi' : \dN \to \dN$ by $\phi'(t w) =
\phi(t) \overline{\phi}(w)$ one easily checks that $\phi' \in
\Aut(\dN,\{\dN_\sigma\})$. Thus it suffices to prove that $\psi =
\phi \phi'^{-1} \in \Aut(\dN,\{\dN_\sigma\})$. Since $\psi$ is
the identity on $\dT$ it follows that $f(w) = \psi(w) w^{-1}$ defines
an element in $\Der(W,\dT)$. Now $H^1(\left<\sigma\right>;\dT) = 0$ so
$f(\sigma) \in (1-\sigma)(\dT) = T^-_0(\sigma)$ and hence $\psi \in
\Aut(\dN,\{\dN_\sigma\})$ as desired.
\end{proof}

\begin{proof}[Proof of Theorem~\ref{mainthm2cg}]
First consider the case $p=2$. The proof of the first part of
Theorem~\ref{mainthmlie} carries over verbatim to give a proof of the
first part of Theorem~\ref{mainthm2cg}. So to finish the proof in this
case we just need to see that $\Phi: \pi_0(\Aut(BX)) \to \Out(\dN)$
has image contained in $\Out(\dN,\{\dN_\sigma\})$, which basically
follows from the definitions:
For any homotopy equivalence $\phi: BX \to BX$ there exists by
\cite[Lem.~4.1 and Prop.~5.1]{AGMV03} an automorphism $\breve \phi: \dN
\to \dN$, making the diagram
$$
\xymatrix{
B\dN \ar[r]^-{B\breve\phi} \ar[d]^-i & B\dN \ar[d]^-i \\
BX \ar[r]^-\phi & BX
}
$$
homotopy commute, and $\breve\phi$ depends only on
the homotopy class of $\phi$, up an inner automorphism of $\dN$. Let
$\dN(\sigma) = C_{\dN}(\dT^+_0(\sigma))$ and $BX(\sigma) =
B\cC_X(\dT^+_0(\sigma)) = \map(B\dT_0^+(\sigma),BX)_i$, where
$i:B\dT_0^+(\sigma) \to BX$ is the inclusion. By definition $\breve\phi(\dN(\sigma)) = \dN({\breve\phi
\sigma \breve\phi^{-1}})$ and the diagram
$$
\xymatrix{
B\dN(\sigma) \ar[d]^{B{\breve \phi}} \ar[r] & BX(\sigma)
\ar[d]_{f\mapsto \phi \circ f\circ B\breve\phi^{-1}} \ar[r] &  BX
\ar[d]^\phi \\
B\dN(\breve\phi \sigma \breve\phi^{-1}) \ar[r]& BX(\breve\phi \sigma \breve\phi^{-1}) \ar[r] & BX
}
$$ 
commutes up to homotopy. This, together with
the definition of $\dN_\sigma$ (see Section~\ref{normalext-sec} and
\cite[Def.~9.5]{DW05}) makes it clear that $\breve \phi(\dN_\sigma) =
\dN_{\breve\phi \sigma \breve\phi^{-1}}$, which is what we needed to prove.

Next consider the case $p$ odd. In this case $\dnu(\D) = \dN = \dT \rtimes
W$ and, e.g.\ by \cite[Prop.~5.2]{AGMV03} there is an exact sequence
$1 \to H^1(W;\dT) \to \Out(\dN) \to N_{\Aut(\dT)}(W)/W \to 1$. Since the
reflections have order prime to $p$ the coroots are given by
$\Z_p b_\sigma = \im(L \xrightarrow{1-\sigma} L)$, so $\Aut(\D) =
N_{\Aut(\dT)}(W)$. This proves the existence of the short exact
sequence in the theorem. By \cite[Thm.~3.3]{kksa:thesis} we have
$H^1(W;\dT) = 0$ so the sequence is canonically split. The theorem now
follows from Lemma~\ref{autlemmapodd}.
\end{proof}

\begin{rem} \label{extclassrem}
It is easy to see that for a $\Z$-root datum $\D$ one has an exact sequence
$$
1 \to \Der(W,T) \to \Aut(\nu(\D)) \to {}^\gamma N_{\Aut(T)}(W) \to 1,
$$
where ${}^\gamma N_{\Aut(T)}(W)$ denotes the subgroup of
$N_{\Aut(T)}(W)$ which fixes the extension class $\gamma$ of the extension
$1 \to T \to \nu(\D) \to W \to 1$ (cf.\ \cite[\S 5]{AGMV03} for
details). Hence Theorem~\ref{mainthmlie} identifies
$\Aut(\D)$ with ${}^\gamma N_{\Aut(T)}(W)$ as subgroups of
$\Aut(T)$. Similarly Theorem~\ref{mainthm2cg} shows that, for a
$\Z_p$-root datum $\D$, $\Aut(\D)$ equals ${}^\gamma N_{\Aut(\dT)}(W)$
as subgroups of $\Aut(\dT)$.
\end{rem}


\section{Construction of $B\textnormal{aut}(\D)$} \label{kinv-sec1}

The goal of this section is, for a $\Z_p$-root datum $\D$, to
introduce the space $B\aut(\D)$ and a refined Adams-Mahmud map
$\Phi: B\Aut(BX) \to B\aut(\D)$ modifying the Adams-Mahmud map
$B\Aut(BX) \to B\Aut(B\N_X)$ constructed in \cite[Lem.~3.1]{AGMV03}.

As mentioned in the introduction, if $p$ is odd, we just  set $B\aut(\D) =
B\Aut(B\N_\D)$, where $B\N_\D = (B^2L)_{hW}$, and let $\Phi: B\Aut(BX) \to B\Aut(B\N_X) \simeq
B\aut(\D)$ be the standard Adams-Mahmud map. 
For $p$ odd we also define $B\daut(\D) = B\Aut(B\dnu(\D))$ so that
we have a canonical map $B\daut(\D) \to B\aut(\D)$ which is a partial
$\F_p$-completion, leaving
the fundamental group unchanged, cf.\ \cite[Ch.~VII, 6.8]{bk}. We now
embark on giving the definitions for $p=2$.

For any $p$, define the {\em discrete center} of $\D$ as 
$$
\dZ (\D) = \bigcap_{\sigma} S(\sigma),
$$
where $S(\sigma)$ is the singular set corresponding to $\sigma$, cf.\
\eqref{singset}. In particular 
\cite[Thm.~7.5]{dw:center} shows that $\dZ(X) \cong \dZ(\D_X)$ for a
connected $p$-compact group $X$. The space $B^2\cZ(\D)$ is defined as
the $\F_p$-completion of $B^2\dZ(\D)$, and the {\em center} $\cZ(\D)$
is defined as the double loop space of this space.

\begin{lemma} \label{centercorrection}
If a reflection $\sigma$ in a $\Z_p$-root datum $\D = (W,L,\{\Z_p
b_\sigma\})$ satisfies $S(\sigma)\neq \dT^+(\sigma)$ then $p=2$ and
$\sigma$ is a reflection in a direct
factor of $\D$ isomorphic to $\D_{\SO(2n+1)} \otimes_\Z \Z_2$.
Furthermore
$$
Z(\dnu(\D)) = \dT^W = \dZ(\D) \oplus V
$$
as $\Out(\D)$- and as $\Out(\dnu(\D))$-modules, where $V=0$ unless
$p=2$ in which case $V \cong (\Z/2)^s$ where $s$ is the number of direct
factors of $\D$ isomorphic to $\D_{\SO(2n+1)} \otimes_\Z \Z_2$.
\end{lemma}

\begin{proof}
If $p$ is odd, then by \cite[Rem.~7.7]{dw:center} $S(\sigma) = \dT^+(\sigma)$ and
by \cite[Thm.~7.5]{dw:center} $Z(\dnu(\D)) = \dT^W = \dZ(\D)$. Hence we can
assume that $p=2$. Write $\D = \D_1 \times \D_2$, where
$\D_2 = (W_2, L_2, \{\Z_2 b_\sigma\})$ is a direct product of copies of
$\D_{\SO(2n+1)\twocom}$,  and $\D_1 = (W_1, L_1,\{\Z_2 b_\sigma \})$
does not contain any direct factors isomorphic to  $\D_{\SO(2n+1)\twocom}$. 

By \cite[Prop.~7.4]{DW05} or \cite[Thm.~11.1]{AGMV03}, we can
write $\D_1$ as a direct product of a $\Z_2$-root datum of Coxeter
type and a number of copies of $\D_{\DI(4)}$. Since $\dT^+_0(\sigma) =
S(\sigma) = \dT^+(\sigma)$ for any reflection $\sigma\in W_{\DI(4)}$
and $Z(\dnu(\D)) = \dT^{W_{\DI(4)}} = 0$ in this case, we can
furthermore assume that $\D_1$ is of Coxeter type. Hence, by
\cite[Prop.~3.1(ii)]{matthey02} (see also
\cite[Prop.~3.2(vi)]{JMO95} or \cite[\S 4]{osse97}), $S(\sigma) =
T_0^+(\sigma)$ for all reflections $\sigma \in W_1$. In particular
$\dT_1^{W_1} = \dZ(\D_1)$ whereas $\dZ(\D_2) = 0$ and $\dT_2^{W_2}
\cong (\Z/2)^s$, where $s$ is the number of direct factors in $\D_2$.

By combining Remark~\ref{extclassrem} with \cite[Prop.~5.4]{AGMV03} we
see that $\Out(\D) = \Out(\D_1) \times \Out(\D_2)$ so $\dT^W =
\dT_1^{W_1} \times \dT_2^{W_2} = \dZ(\D) \times V$ splits as an
$\Out(\D)$-module in the manner indicated in the lemma.
Since the $\Out(\dnu(\D))$-action factors through the
$\Out(\D)$-action, the above splitting is also a splitting as
$\Out(\dnu(\D))$-modules.
\end{proof}

From now on, let $\D$ be a $\Z_2$-root datum and let $\dN = \dnu(\D)$ be the
associated discrete maximal torus normalizer. Let $B\N = B\N_\D$ denote the
fiber-wise $\F_2$-completion of $B\dN$ with respect to the fibration $B\dN
\to BW$ \cite[Ch.~I, \S 8]{bk}.
Define $\breve Y$ to be the covering space of $B\Aut(B\dN)$
corresponding to the subgroup $\Out(\dN,\{ \dN_\sigma\})$ of
$\Out(\dN) = \pi_1(B\Aut(B\dN))$.

By Lemma~\ref{centercorrection}, $\pi_2(\breve Y) \cong
\dZ(\D) \oplus V$ as $\pi_1(\breve Y)$-modules, since $\pi_1(\breve Y) =
\Out(\dN,\{ \dN_\sigma\}) \xrightarrow{\cong} \Out(\D)$ by
Theorem~\ref{mainthm2cg}. Define $B\daut(\D)$ to be the space obtained
from $\breve Y$ by first attaching $3$-cells to kill exactly $V
\subseteq \pi_2(\breve Y)$, which we can do since $V$ is a
$\pi_1(\breve Y)$-invariant subgroup of $\pi_2(\breve Y)$, and then taking the second
Postnikov section, killing all homotopy groups in dimensions three and
higher. By construction $B\daut(\D)$ satisfies
$$
\pi_i(B\daut(\D)) =
\begin{cases}
\Out(\D) & \text{for $i=1$,} \\
\dZ(\D) & \text{for $i=2$,} \\
0 & \text{otherwise.}
\end{cases}
$$ 
We define $Y$ to be the covering space of $B\Aut(B\N)$ with respect
to the subgroup $\Out(\dN,\{ \dN_\sigma\})$ of $\pi_1(B\Aut(B\N))
\xleftarrow{\cong} \Out(\dN)$, so that we have a canonical map $\breve Y \to
Y$. Define $B\aut(\D)$ to be the partial $\F_2$-completion
of $B\daut(\D)$ leaving the fundamental group unchanged. By
construction $B\aut(\D)$ has universal cover homotopy
equivalent to $B^2\cZ(\D)$ and in particular $B\aut(\D)$ only
has homotopy groups in dimensions $1$, $2$, and $3$.

Since the maps $\breve Y \to Y$ and $B\daut(\D) \to B\aut(\D)$ are
partial $\F_2$-completions, the map $\breve Y \to B\daut(\D)$ induces a
map $Y \to B\aut(\D)$, well defined up to homotopy, by the universal
property of the partial $\F_2$-completion. On homotopy groups this map just kills the summand $V$ in $\pi_2(Y)$.

By Theorem~\ref{mainthm2cg} the image of the map $\Phi:
\pi_1(B\Aut(BX)) \to \pi_1(B\Aut(B\N))$ is contained in $\Out(\dN,\{\dN_\sigma\})$, so it lifts uniquely
to a basepoint-preserving map $B\Aut(BX) \to Y$. We define
$$
\Phi: B\Aut(BX) \to B\aut(\D)
$$
to be the composite $B\Aut(BX) \to Y \to B\aut(\D)$.

\medskip

We note the following result, which we will need in the proof of
Theorem~\ref{modelthm}.

\begin{prop} \label{Bautdproduct}
If $\D = \D_1 \times \D_2$ is a product of $\Z_p$-root data such that
$\D_1$ and $\D_2$ have no direct factors in common, then we have a
canonical map
$B\aut(\D_1) \times B\aut(\D_2) \to B\aut(\D)$ which is a homotopy equivalence.
\end{prop}

\begin{proof}
There is a natural map $B\Aut(B\N_{\D_1}) \times B\Aut(B\N_{\D_2}) \to
B\Aut(B\N_\D)$ which by \cite[Prop.~5.4]{AGMV03} is a homotopy equivalence for $p$ odd.
Assume that $p=2$. Here, by Remark~\ref{extclassrem} and
\cite[Prop.~5.4]{AGMV03}, $\Out(\D_1) \times \Out(\D_2)
\xrightarrow{\cong} \Out(\D)$, and hence we have an isomorphism
\begin{equation} \label{natural}
\Out(\dnu(\D_1),\{\dnu(\D_1)_\sigma\}) \times
\Out(\dnu(\D_2),\{\dnu(\D_2)_\sigma\}) \xrightarrow{\cong}
\Out(\dnu(\D),\{\dnu(\D)_\sigma\})
\end{equation}
Passing to covers of the natural map above with respect to the isomorphic
subgroups of the fundamental groups given by the isomorphism \eqref{natural}
produces a canonical map $Y_{\D_1} \times Y_{\D_2} \to Y_\D$ which by
construction is a homotopy equivalence, since it induces an isomorphism on
all homotopy groups. Killing the $\Z/2$ summands in $\pi_2$
corresponding to direct factors isomorphic to $\D_{\SO(2n+1)\twocom}$ produces a
homotopy equivalence $B\aut(\D_1) \times B\aut(\D_2) \to B\aut(\D)$ as
wanted.
\end{proof}


\section{Proof of Theorems \ref{kinvthm} and \ref{modelthm}} \label{kinv-sec2}

Theorems~\ref{kinvthm} and \ref{modelthm} will follow
easily from the results of the previous section, combined with the
following theorem, which is an analog of a classical theorem for
Lie groups by de Siebenthal \cite[Ch.~I, \S 2, no.\ 2]{desiebenthal56}.

\begin{thm} \label{plingeling}
\InsertTheoremBreak
\begin{enumerate}
\item \label{pling}
Let $\D$ be a $\Z$-root datum and let $N = \nu(\D)$. Then the exact sequence
$$
1 \to \Inn(N) \to \Aut(N,\{N_\sigma\}) \to \Out(N,\{N_\sigma\}) \to 1
$$
is split.
\item \label{pling2}
Let $\D$ be a $\Z_2$-root datum of Coxeter type. Then the same statement
holds if we replace $N=\nu(\D)$ by $\dN=\dnu(\D)$.
\item \label{plingodd}
Let $\D$ be a $\Z_p$-root datum of Coxeter type, $p$ odd, with
associated discrete maximal torus normalizer $\dN = \dnu(\D)$. Then
the exact sequence
$$
1 \to \Inn(\dN) \to \Aut(\dN) \to \Out(\dN) \to 1
$$
is split.
\end{enumerate}
\end{thm}

\begin{proof}
\eqref{pling} We define a splitting $s: \Out(N,\{N_\sigma\}) \to
\Aut(N,\{N_\sigma\})$ as follows. First fix a simple system of roots $B$ and
let $S$ be the corresponding set of simple reflections. We also fix a
set of elements $\{x_\sigma\}_{\sigma\in S}$ with $x_\sigma \in
N_\sigma\setminus T_0^-(\sigma)$. Let $[\phi]\in
\Out(N,\{N_\sigma\})$ be an element represented by an automorphism
$\phi \in \Aut(N,\{N_\sigma\})$. Since the action of $W$ on the set
of simple systems of roots is simply transitive (e.g. see
\cite[Thm.~1.8]{humphreys90}), we can find a unique element $w\in W$
with $\phi(B) = w(B)$. By composing $\phi$ with an appropriate
inner automorphism of $N$ we thus see that there exists a representative
$\phi'$ of $[\phi]$ with $\phi'(B) = B$. Moreover
$\phi'$ is unique up to multiplication by an inner automorphism
given by conjugation by an element in $T$.

Since $\phi'(S) = S$ and $\phi' \in \Aut(N,\{N_\sigma\})$ we
have $\phi'(x_\sigma) \in N_{\phi'(\sigma)}\setminus
T_0^-(\phi'(\sigma))$. Hence $\phi'(x_\sigma) = t_{\phi'(\sigma)}\cdot
x_{\phi'(\sigma)}$ for certain elements $t_\sigma \in
T_0^-(\sigma)$, $\sigma\in S$, depending on $\phi'$. Let $\phi'' = c_t \circ
\phi'$, where $c_t\in \Inn(N)$ denotes conjugation by an element
$t\in T$. We now have
$$
\phi''(x_\sigma) = t\cdot t_{\phi'(\sigma)}\cdot
x_{\phi'(\sigma)}\cdot t^{-1} = t_{\phi'(\sigma)}\cdot t\cdot
\left(\phi'(\sigma)(t)\right)^{-1}\cdot x_{\phi'(\sigma)}.
$$
Since $t_{\phi'(\sigma)} \in T_0^-(\phi'(\sigma))$ for
$\sigma\in S$ it now follows from
Lemma~\ref{Tsurjlemma} that we can find $t\in T$ such that
$\phi''(x_\sigma) = x_{\phi'(\sigma)}$. Moreover it is clear
that such a $t\in T$ is unique up to multiplication by an element in
$T^W = Z(N)$. We conclude that any element $[\phi]\in
\Out(N,\{N_\sigma\})$ has a unique representative $\phi''\in
\Aut(N,\{N_\sigma\})$ with $\phi''(B) = B$ and $\phi''(x_\sigma) =
x_{\phi''(\sigma)}$. Hence the assignment
$s([\phi]) = \phi''$ will be a group homomorphism and define the
desired splitting.

\eqref{pling2} and \eqref{plingodd}: When $\D$ is a $\Z_p$-root
datum of Coxeter type the proof above goes through verbatim to show
that the short exact sequence $1 \to \Inn(\dN) \to \Aut(\dN,\{\dN_\sigma\}) \to
\Out(\dN,\{\dN_\sigma\}) \to 1$ is split. This proves the result since
$\Aut(\dN,\{\dN_\sigma\}) = \Aut(\dN)$ for $p$ odd by
Lemma~\ref{autlemmapodd}.
\end{proof}

\begin{proof}[Proof of Theorem~\ref{modelthm}]
Write $\D = \D_X = (W, L, \{\Z_p b_\sigma\})$. We constructed the map
$\Phi: B\Aut(BX) \to B\aut(\D)$ in Section~\ref{kinv-sec1}, and from
the construction it is clear that $\Phi$ induces an isomorphism on
$\pi_i$ for $i>1$, and that $B\aut(\D)$ is the total space of a
fibration $B^2\cZ(\D) \to B\aut(\D) \to B\Out(\D)$. The remaining
claim is that this fibration is split, which we now prove. The space
$B\aut(\D)$ is by definition the partial $\F_p$-completion of the
total space of the fibration
\begin{equation} \label{BdautDfib}
B^2\dZ(\D) \to B\daut(\D) \to B\Out(\D),
\end{equation}
and it is hence enough to see that this fibration splits. By
\cite[Thm.~11.1]{AGMV03} (cf.\ also \cite[Prop.~1.12]{DW05}) we may
write $\D = \D_1 \times \D_2$, where $\D_1$ is of Coxeter type and
$\D_2$ is a direct product of exotic $\Z_p$-root data. Since $\D_1$ and $\D_2$
have no direct factors in common, it follows from
Proposition~\ref{Bautdproduct} combined with Remark~\ref{extclassrem}
and \cite[Prop.~5.4]{AGMV03} that the fibration \eqref{BdautDfib} is the
product of the corresponding fibrations for $\D_1$ and $\D_2$. By
\cite[Thm.~11.1]{AGMV03} we have $\dZ(\D_2) = 0$, so it suffices to
prove that \eqref{BdautDfib} splits when $\D$ is of Coxeter type.

To prove this consider the composite
$$
B\Out(\dN,\{\dN_\sigma\}) \to B\Aut(\dN,\{\dN_\sigma\}) \to
{\breve Y} \to B\daut(\D).
$$
Here the first map is induced by the splitting constructed in
Theorem~\ref{plingeling}, the second map is the obvious
map sending an automorphism to the induced self-homotopy equivalence
of classifying spaces which factors through ${\breve Y}$,  and
the third map is the map constructed in Section~\ref{kinv-sec1}.
By construction this is a section to the map $B\daut(\D) \to
B\Out(\D) \simeq B\Out(\dN,\{\dN_\sigma\})$ as desired.
\end{proof}

\begin{proof}[Proof of Theorem~\ref{kinvthm}]
We continue with the notation of the previous proof. By construction
the square
\begin{equation} \label{pbs}
\xymatrix{B\Aut(BX) \ar[r] \ar[d] & B\aut(\D) \ar[d] \\
B\Out(BX) \ar[r] & B\Out(\D)}
\end{equation}
is a homotopy pullback square, so by the universal property of the
homotopy pullback, the fibration $B^2\cZ(X) \to B\Aut(BX) \to
B\Out(BX)$ has a section since $B^2\cZ(\D) \to B\aut(\D) \to
B\Out(\D)$ has one.

The fact that any homomorphism $\Gamma \to \Out(BX)$ lifts to an action
follows directly from the splitting. To see the parametrization in the case
where $\widetilde{H}^*(\Gamma;\Q) =0$, 
note that, by the above pullback square \eqref{pbs}, liftings to
$B\Aut(BX)$ agree with lifts
$$
\xymatrix{
 & B\aut(\D) \ar[d] \\
B\Gamma \ar[r] \ar@.[ur] \ar@.[ur] & B\Out(\D).
}
$$
But since $B\Gamma$ is assumed to be rationally trivial, the fibration sequence
$$
B^2(L^W \otimes \Q) \to B\daut(\D) \to B\aut(\D)
$$
shows that these lifts
agree with lifts in the same diagram, but with $B\aut(\D)$ replaced by
$B\daut(\D)$. Because \eqref{BdautDfib} is split, obstruction theory
now directly implies that the set of lifts is a non-trivial
$H^2(\Gamma;\dZ(\D))$-torsor.
\end{proof}

\begin{rem}
Theorem~\ref{kinvthm} in particular implies that the second $k$-invariant of
$B\Aut(BX)$ vanishes, since $P_2(B\Aut(BX)) \to P_1(B\Aut(BX))$ has a
section, where $P_n(\cdot)$ denotes the $n$th Postnikov piece.  The
only possible non-zero $k$-invariant of $B\Aut(BX)$ is therefore the
third, which however need not vanish. This can occur for all $p$ (as
well as for compact Lie groups). The third $k$-invariant vanishes if
and only if $P_{3}B\Aut(BX) \to P_{2}B\Aut(BX)$ has a section, which
is equivalent to $B\cZ(X)$ being a product of Eilenberg-Mac\,Lane spaces
{\em as an $\Out(BX)$-space}, which need not be the case. (By
$\Out(G)$-spaces we mean the model category of spaces with an
$\Out(G)$-action, where the weak equivalences are the
$\Out(G)$-equivariant maps which are non-equivariant homotopy
equivalences; see e.g, \cite[VI.4]{GJ99}.)
The point is that while $\dZ(X) \to \dZ(X)/C$ always has a section,
where $C$ is the largest divisible subgroup of $\dZ(X)$, it need {\em
  not} have an $\Out(BX)$-equivariant section.

A concrete example is given by 
$G = (S^1 \times \SU(n) \times \SU(n))/\Delta$, where $\Delta$ is the
``diagonal'' central subgroup generated by $(\zeta, \zeta
I,\zeta I)$, $\zeta=e^{2\pi i/n}$---we leave the details to the
reader.
\end{rem}


\bibliographystyle{plain}
\bibliography{../pcg/poddclassification}
\end{document}